\documentclass[acmsmall]{article}
\usepackage[a4paper, total={6in, 8in}]{geometry}
\usepackage{amsmath,amssymb,amsthm}
\usepackage{hyperref}
\usepackage{graphicx}
\usepackage{color}
\usepackage{algorithm}
\usepackage{algpseudocode}
\usepackage{multirow}
\usepackage{authblk}
\usepackage{caption}
\usepackage{subcaption}

\newtheorem{lemma}{Lemma}
\newtheorem{theorem}{Theorem}
\newtheorem{proposition}{Proposition}
\newtheorem{corollary}{Corollary}
\newtheorem{remark}{Remark}

\newcommand{\abs}[1]{\lvert #1 \rvert}
\newcommand{\brac}[1]{\left(#1\right)}
\newcommand{\sbrac}[1]{\left[#1\right]}

\newcommand{\Nats}{\mathbb{N}}
\newcommand{\Ints}{\mathbb{Z}}
\newcommand{\Reals}{\mathbb{R}}

\newcommand{\prob}{\mathbb{P}}
\newcommand{\E}{\mathbb{E}}

\begin{document}

\title{A Model of Job Parallelism for Latency Reduction
in Large-Scale Systems}

\author[$*$]{A. Ganesh}
\author[$\dag$]{A. Mukhopadhyay}
\affil[$*$]{School of Mathematics, University of Bristol}
\affil[$\dag$]{Department of Computer Science, University of Warwick}

\date{}

\maketitle

\begin{abstract}
  Processing computation-intensive jobs at multiple processing cores in parallel is essential in many real-world applications. In this paper, we consider an idealised model for job parallelism in which a job can be served simultaneously by $d$ distinct servers. The job is considered complete when the total amount of work done on it by the $d$ servers equals its size. We study the effect of parallelism on the average delay of jobs. Specifically, we analyze a system consisting of $n$ parallel processor sharing servers in which jobs arrive according to a Poisson process of rate $n \lambda$ ($\lambda <1$) and each job brings an exponentially distributed amount of work with unit mean. Upon arrival, a job selects $d$ servers uniformly at random and joins all the chosen servers simultaneously. We show by a mean-field analysis that, for fixed $d \geq 2$ and large $n$, the average occupancy of servers is $O(\log (1/(1-\lambda)))$ as $\lambda \to 1$ in comparison to $O(1/(1-\lambda))$ average occupancy for $d=1$. Thus, we obtain an exponential reduction in the response time of jobs through parallelism. We make significant progress towards rigorously justifying the mean-field analysis.
\end{abstract}




\section{Introduction} \label{sec:intro}

Consider a model consisting of $n$ parallel servers,
each having one unit of processing capacity and employing the processor-sharing (PS) scheduling discipline. Jobs arrive according to a Poisson process with rate $n \lambda$ and each job brings a random amount of work, exponentially distributed with unit mean. Upon arrival of a job, $d$ servers are chosen uniformly at random.
We assume an idealised model where these $d$ servers process the job simultaneously and the job leaves the
system as soon as the total amount of work done on it
by all $d$ servers equals its size. Since the servers
operate according to the PS discipline,
the instantaneous rate obtained by a job at a given server
is a fraction of the server's total capacity, which is shared equally by all the jobs present at that server. Our goal is to study the average delay of jobs in the system as a function of the parameter $d$, which captures the degree of parallelism in the system.
We focus on the limiting regime where $n$ tends to infinity. We note that for $d=1$ the system
reduces to $n$ independent parallel servers for which
the average delay is known to be $1/(1-\lambda)$. However, for $d \geq 2$, the servers' state are not independent of each other, since at any given instant a server interacts with the subset of servers with which it shares jobs. This makes the analysis for  $d \geq 2$ more challenging.

The motivation to study this idealised model comes from modern 
data centres, which often handle computing jobs that are highly parallelisable. For example, machine learning jobs such as 
TensorFlow~\cite{tensorflow_2016} and scientific computing jobs 
often require large numbers of parallel cores for efficient processing~\cite{GoogleBorg_2015}. Other examples of large-scale 
parallelism include coded file retrieval systems~\cite{lee2017} 
and parallel computing systems implementing the MapReduce  framework~\cite{mapreduce}. The average time a job spends
in such a system crucially depends on the degree of parallelism used, and on how the capacities of the processing cores are shared by jobs present in the system.

Various queueing models have been proposed to study
the effect of parallelism and scheduling on job delay.
Such models are, in general, difficult to analyze due to the dependency between servers induced by parallel jobs. 
Many recent works analyse models in which jobs are divided into tasks 
and processed at multiple servers~\cite{joshi_latency_redundant,weina_fj_2018,rizk_2016}; 
a job is considered to be completed when all its tasks finish processing. These are analysed as generalised Fork-Join (FJ) 
queuing models.
Due to the assumption that servers operate a first-come-first-served (FCFS) policy, in these models a job may not be served simultaneously at multiple servers. Thus, the dependence among servers in these models is {\em weak}, in the sense that the servers are only coupled at arrival instants of jobs. Another class of 
models~\cite{afanaseva_2020,rumayanstev_2017,weina_multiserv_2021} 
studies systems consisting of a single central queue and multiple parallel servers. Jobs requiring multiple servers arrive into the central queue. If the required number of servers is available upon arrival, then the job joins all the servers simultaneously; otherwise, it waits in the queue until the required number of servers become available. 
In the models discussed above, servers do not have individual queues and as a result can serve only one job at a time. Our model differs in allowing each server to process multiple jobs simultaneously, and each job to receive service from multiple servers simultaneously.

An important question in this context is how to select $d$, the degree of parallelism in the system. At one extreme, when $d=1$, the average delay is known to be $\frac{1}{1-\lambda}$, whereas the other extreme of $d=n$ is equivalent to perfect resource pooling and has average delay of $\frac{1}{n(1-\lambda)}$. But this comes at the cost of greatly increased communication between servers. To what extent can the benefits of resource pooling be realised with small values of $d$?
Our key finding is that, in the large $n$ limit for any fixed $d \geq 2$, the average delay of jobs scales as $O(\frac{1}{\lambda}\log \frac{1}{1-\lambda})$ as $\lambda$ increases to $1$. This implies that even with a small degree of parallelism it is possible to obtain an exponential reduction in the average delay of jobs in heavy traffic. A small value of $d$ also implies that the communication overhead can be kept low. Our main contributions are the following.

\subsection{Contributions}

In Section~\ref{sec:mean-field}, we present a  mean-field approximation of the steady-state average delay of jobs. To obtain this approximation, we study the evolution of a tagged queue, assuming that all the servers with which the tagged server shares a job are in their equilibrium distribution at all times. This assumption yields a Markovian evolution of the tagged queue with transition rates which depend on the equilibrium distribution. Solving for the equilibrium distribution of this Markov process yields a fixed point equation. We show that this fixed point equation has a unique solution for every $\lambda<1$. We solve it to obtain the equilibrium occupancy distribution and the mean occupancy per server. We show that the occupancy is bounded by a constant multiple of $-\log(1-\lambda)$ as $\lambda$ increases to 1; this is in stark contrast to the $1/(1-\lambda)$ scaling in the $d=1$ case, and is reminiscent of the double exponential decay seen in the supermarket queueing model~\cite{dobrushin96,mitzenmacher01}. Using similar assumptions, we also obtain the evolution equations for the transient distribution of the mean-field model.


In Section~\ref{sec:numerical}, we present evidence from extensive computer simulations to support the conjecture that the solution to the mean-field equations yields the invariant marginal occupancy distribution. We overlay plots of the empirical cumulative distribution function (cdf) from simulations on the theoretical cdf from the mean-field equations, for systems consisting of $n=100$ and $n=500$ servers, and different values of the arrival rate, $\lambda$. The plots show that the empirical and theoretical cdfs are very close to each other. Likewise, the mean value from simulations is very close to the theoretical mean. We also present simulation results which indicate that for large system sizes the stationary queue-length distribution of our model is insensitive to the type of service-time distribution as long as its mean remains unchanged. This is a very desirable property for system designers as it implies that system performance does not change with the change in the service time distribution of arriving tasks.

Next, we turn to the problem of establishing rigorously that the equilibrium occupancy distribution in the system with $n$ servers converges to the mean-field prediction as $n$ tends to infinity. The analysis turns out to be significantly more challenging than other related mean-field models due to the fact
that the dependence among the servers' states is much stronger
in our model. Indeed, in our model, a server's state is coupled with every other server with which it shares a job, for as long as the job is present in the system. Thus, unlike many other mean-field models in the literature, the interaction between servers in our model is not restricted to specific epochs, such as arrival instants. Despite this difficulty, we make significant progress towards rigorously justifying the mean-field model.
In Section~\ref{sec:stability}, we show that the system is stable under the natural condition, $\lambda<1$, that the rate at which work arrives into the system is smaller than the maximum rate at which it can be served. This shows the existence of a unique invariant distribution $\pi^n$ for the number of jobs at each server. 
We further show that the mean number of jobs at a server in equilibrium is bounded uniformly in $n$. This implies that the sequence of distributions $(\pi^n)_{n\in \Nats}$, is tight, and hence that it has subsequential weak limits. We conjecture that the limit is unique.

Next, in Section~\ref{sec:monotone}, we show that the system, started empty, converges monotonically (in the stochastic order induced by the componentwise partial order) to its invariant distribution.
This result is essential in studying the convergence of mean-field
models to their corresponding fixed point.
Using this monotonicity result, in Section~\ref{sec:unif_conv} we show that the speed of convergence of the marginal occupancy distribution to $\pi^n$ is uniform in $n$. This result is a significant step towards showing propagation of chaos, i.e., that any fixed number of queues become asymptotically independent as the system size $n$ tends to infinity. It shows that
the marginal occupancy distribution for large $n$ 
and large $t$ can be close both
to the stationary distribution and the distribution of a cavity process. In Section~\ref{sec:discussion}, we discuss
in more detail how our theoretical results could be used
to prove the convergence to a cavity process.
However, it remains an open problem to define the cavity process and to show that for every fixed $t$ the system converges to this cavity process as $n \to \infty$.

Our final contribution, discussed in Section~\ref{sec:static}, is to study a static version of the model, in which $\lambda n$ unit-sized jobs are each replicated on a set of $d$ servers out of $n$, chosen independently and uniformly at random. Here, $\lambda>0$ is a fixed constant which does not depend on $n$. We consider the problem of scheduling the servers so as to minimise the makespan, namely, the time until all jobs are completed. We show that with high probability, i.e., with probability tending to 1 as $n$ tends to infinity, the makespan is bounded uniformly in $n$, as a function of $\lambda$. In other words, the jobs can be scheduled in such a way that the time for even the most heavily loaded server to complete its tasks does not grow without bound. This result prompts the question of whether one can find better algorithms in the dynamic setting than that of sharing server capacity equally amongst all jobs present at a server. This is an open problem.

\subsection{Related work}

Load-balancing in large systems has attracted a great deal of research. Here, we selectively discuss works which are related to the model or the analysis techniques studied in this paper.

In the context of load balancing in parallel server systems, the SQ($d$) algorithm was first analyzed independently in~\cite{dobrushin96, mitzenmacher01} as an alternative to the classical Join-the-Shortest-Queue (JSQ) scheme to reduce the overhead of job dispatching. In the SQ($d$) algorithm, an incoming job is sent to the shortest among $d$ queues chosen uniformly at random. Using mean-field techniques it was shown that in the limit as the number of servers tends to infinity, the marginal queue length distribution exhibits double exponential decay in the tail for $d \geq 2$; this is in contrast to assigning new jobs uniformly at random, which gives rise to independent $M/M/1$ queues, which have a geometric queue length distribution (i.e., the queue length decays exponentially). A closely related algorithm, denoted LL($d$), in which jobs join the least loaded of $d$ servers was analyzed in the heavy traffic regime in~\cite{leastloaded2018}; this algorithm requires knowledge of the workload associated with each job. The analysis was extended to a number of other load-balancing policies in~\cite{least_loaded_heavy_traffic}, in the heavy traffic regime.
Bramson {\em et al.}~\cite{bramson2012} analyzed the SQ($d$) and LL($d$) algorithm for general service time distributions using the cavity method from statistical physics. They established that in the mean-field limit any finite number of queues become independent of each other. The approach adopted in this paper closely follows the general framework introduced in their work.

A generalisation of SQ($d$) to multi-component jobs is studied in~\cite{shneer2021}. Here, jobs consist of $k$ components; $d$ servers are sampled, and components are assigned to the least-loaded among these servers so as to balance their loads. This combines elements of parallelism and load-balancing.


The above approaches assume knowledge of server loads when assigning jobs. An alternative that has been extensively studied recently involves placing replicas of a job at multiple servers, in the hope that at least one of them will have low load.
Suppose each arriving job is replicated at $d$ servers chosen uniformly at random. When one of these servers either starts working on a replica, or completes working on it, then replicas of the job at all $d-1$ other servers are immediately removed; the former policy is called cancel-on-start (CoS) and the latter cancel-on-complete (CoC). These policies have been studied in~\cite{gardner16,gardner17,ayesta18}, which give exact closed-form expressions for the invariant distribution and mean response time in small systems. These expressions involve a state space representation whose complexity grows rapidly in the number of servers and it is not easy to obtain insights into large-system asymptotics. Stability of replication policies was studied in~\cite{stability_redundancy_fcfs} for the FCFS service discipline, and in~\cite{anton21} for several different service disciplines, including the one considered in this paper. Both stability and tail behaviour of such policies under the PS service discipline are analyzed in~\cite{tail_redundancy_ps}. An alternative approach to dealing with server loads which are not known in advance is studied in~\cite{ganesh2012}. Here, jobs are assigned to a random server upon arrival, but periodically sample other servers, and migrate if the sampled server is serving fewer jobs. The paper presents a mean-field analysis of the resulting model, as well as a rigorous justification of it.

A variant of replication, known as coded computation, exploits the fact that many computationally intensive jobs are intrinsically parallelisable, i.e., they can be split into a large number of tasks which can be executed in parallel. Moreover, redundancy can be added via coding so that the job is complete when sufficiently many tasks have been completed. Concretely, a job may be split into $k$ tasks, and a further $m-k$ tasks created via coding. The resulting $m$ tasks can be executed in parallel, and the job finishes as soon as any $k$ tasks are completed. The latency in such a system is studied in~\cite{lee2017, mallick2020}.

\section{System Model and Notation} 
\label{sec:model}

The system consists of $n$ servers, each using processor sharing (PS) as the service discipline and each processing work at unit rate.
This implies that if there are $k$ ongoing jobs at a given server,
then the instantaneous rate at which the server processes each job is $1/k$.
Jobs arrive into the system according to a Poisson processes of rate $n\lambda$. Each job brings in an amount of work which is exponentially distributed with unit mean, independent of all other jobs and of the arrival process. Upon arrival, a job is assigned to $d$ servers, chosen uniformly at random and independent of the past and it is processed in parallel by all the servers to which it is assigned. 
When the total work done on a job by all servers to which it was assigned is equal to the work brought in by that job, the job leaves the system. 
This is equivalent to each job having $d$ copies, each of which is sent to a different server. The servers work on the individual copies independently until the total amount of work done on the job by all servers (this is the sum of the work done on individual copies) equals the job's size.

We shall assume that $\lambda<1$. This condition is necessary for the stability of the system as it ensures that the rate at which work enters the system is strictly smaller than the maximum rate at which work can be processed by the system. We shall later establish that this condition is also sufficient for stability.

We refer to jobs that are processed by a fixed subset of $d$ servers as a {\em class}. Clearly, there are $m=\binom{n}{d}$ classes of jobs. Let $X_i^n(t)$ denote the number of jobs of class $i \in [m]$ in the system at time $t \geq 0$, where we write $[m]$ for the set $\{ 1,\ldots,m \}$. It is easy to see that $X^n()=(X_i^n(),i \in [m])$ is a Markov process on the state space $S^n=\Ints_{+}^{m}$. 
For each class $i \in [m]$, we denote by $\partial_i$
the set of $d$ servers that serves class $i$ jobs and, for each
server $j \in [n]$, we denote by $\partial_j$ the $\binom{n-1}{d-1}$ classes that are served by server $j$. 
Let $Q^n_j(t)=\sum_{i \in \partial_j} X_i^n(t)$ denote the number of ongoing jobs present at server $j$ at time $t$ and let $Q^n(t)=(Q^n_j(t),j \in [n])$.
For simplicity, we refer to $Q_j^n(t)$ as the {\em queue length} of server $j$ at time $t$ although there is no queue at the server and all jobs are processed in parallel. 

We denote by $\pi^{n,0}(t)$ the queue-length distribution of a tagged server at time $t$ starting from the empty state $X^n(0)=\boldsymbol{0}$\footnote{Since the system remains exchangeable starting from the empty state, the queue-length distribution $\pi^{n,0}(t)$ does not depend on the index of the tagged server.} and let $\pi^n$ denote the stationary queue-length distribution of a tagged server (when the system is stable). To analyse the mean response time of jobs for large system sizes, we need to characterise the limit of $\pi^n$ as $n$ tends to infinity and this is the main goal of the paper. The standard technique for establishing 
the limit of $\pi^n$ involves analysing the empirical queue-length process $y^n()=(y^n_k(),k \in \Ints_+)$ defined as $$y^n_k(t)=\frac{1}{n}\sum_{j \in [n]} 1(Q^n_j(t)=k)$$ for each time $t\geq 0$. Typically, the limit of $\pi^n$ coincides with
 $\lim_{t\to\infty}\lim_{n\to\infty} y^n(t)$, which is referred to as the {\em mean-field limit} of $\pi^n$.
However, it is important to note that for the system described above the {\em queue-length process} $Q^n()$ is {\em non-Markovian}, which implies that the empirical measure process $y^n()$ is also non-Markovian. This makes analysing the dynamics of $y^n()$ (and therefore computing $\lim_{t\to\infty}\lim_{n\to\infty} y^n(t)$) significantly more challenging than the  analyses in~\cite{dobrushin96,mitzenmacher01,bramson2012}, as the standard results for density-dependent Markov processes~\cite{kurtz_1970} do not apply.
Nevertheless, in the next section, we present a mean-field limit of $\pi^n$ based on a heuristic framework known as the {\em cavity method} in statistical physics, and analytically characterise the mean-field limit. 



\section{Mean-field analysis} \label{sec:mean-field}

As noted earlier, the queue-length process $Q^n()$ and the empirical queue-length process $y^n()$ are non-Markovian for our system. This makes characterising the mean-field limit of $\pi^n$ difficult using standard techniques. In this section, to characterise the limit $\pi=\lim_{n \to \infty} \pi^n$, we take a heuristic approach based on the cavity method used in statistical physics~\cite{mezard_cavity}.
The cavity method, applied to our model, consists of analysing the invariant distribution of the queue length $Q(t)$ of a single {\em tagged} server in the infinite system under the following assumptions:

\begin{enumerate}
    \item The queues at all servers other than the tagged server are in the steady state (i.e., in the invariant distribution of the process $Q()$) and are independent of the queue length of the tagged server.
    
    \item If a job has a copy present in the tagged server, then the other $d-1$ copies are randomly placed in the rest of the system, i.e., any copy present in the rest of the system has equal probability of being a copy of that job.
\end{enumerate}
Under the above assumptions, the process $Q(t), t \geq 0$ is Markov and has state dependent transition rates which we compute below.

The aggregate arrival rate at any given server in the $n^{\rm th}$ system is $\lambda d$, as the arrival rate into the system is $n\lambda$, and each arrival chooses a subset of $d$ servers uniformly at random; there are $m=\binom{n}{d}$ such subsets, of which $\binom{n-1}{d-1}$ contain the tagged server. Since this is true for all $n$, we take this to be the arrival rate at the tagged queue in the infinite system. Thus, the transition rate of $Q(\cdot)$ from $k$ to $k+1$ is $\lambda d$, for any $k\in \Ints_+$.




Next, for $k\geq 1$, we compute the transition rate from $k$ to $k-1$. Each job in the tagged server is worked on at rate $1/k$; this contributes a total rate of 1 to departures. In addition, each given job is processed at $d-1$ other servers each of which (by the first assumption) has queue length distributed according to the steady state distribution $\pi$ of $Q()$.
Let $Q_{j_1},\ldots,Q_{j_{d-1}}$ denote the stationary random queue lengths seen by
the $d-1$ other copies of a job present at the tagged server.
Then, by the second assumption above, the distributions of $Q_{j_1},\ldots,Q_{j_{d-1}}$ are the same for all jobs whose copy is present at the tagged server. Moreover, since jobs are placed randomly in the system, the probability
that a copy ends up in a queue with length $u$ is proportional to $u \pi_u$ where $\pi_u$ is the $u^{\textrm{th}}$ component of the steady-state distribution $\pi$ of $Q()$. Hence, we have

$$
\prob(Q_{j_l}=u) = \frac{u\pi_u}{\sum_{v=0}^{\infty} v\pi_v} = \frac{u\pi_u}{\E Q}, \quad u\in \Nats, l \in [d-1]
$$
where $\E Q$ denotes the mean queue length in stationarity. 
Since $\prob(Q_{j_l}=0)=0$, it follows that 
\begin{align*}
\E \Bigl( \frac{1}{Q_{j_l}} \Bigr) &= \sum_{u=1}^{\infty} \frac{1}{u} \prob(Q_{j_l}=u) \\
&= \sum_{u=1}^{\infty} \frac{\pi_u}{\E Q} = \frac{1-\pi_0}{\E Q}.
\end{align*}
This is the instantaneous rate at which each copy of a job is worked on by servers other than the tagged server; in other words, it is the service received from the mean field. As this is the case for each job, we conclude that the transition rate from $k$ to $k-1$ is given by $1+(d-1)k\frac{1-\pi_0}{\E Q}$.

We can separately obtain $\pi_0$ by work conservation. The total rate at which work arrives into the system is $n\lambda$, while the total rate at which work is done is $n(1-\pi_0)$, provided the system is stable and a stationary regime exists. Hence $1-\pi_0=\lambda$. 

We now write down transition rates for the queue-length process $Q()$ at the tagged server based on the heuristic arguments presented above. 
Henceforth, by the mean-field model, we refer to the Markov process $Q()$ defined on $\Ints_+$ with the following stationary transition rates: 
\begin{equation}
    \label{eq:trans_rates}
    \begin{aligned}
    &r_{k,k+1}=d\lambda, \quad r_{k+1,k}=1+(k+1)\mu, \quad k\in \Ints_+, \\
    &\mbox{where $\mu=\frac{(d-1)\lambda}{\E Q}$}.
    \end{aligned}
\end{equation}
All other transition rates are zero. We recognise this as an $M/M/1$ queue with reneging, where arrivals occur at rate $d\lambda$, services at rate 1, and each customer reneges independently at rate $\mu$. It is stable for any $\mu>0$, and the invariant distribution solves the local balance equations, $\pi_i q_{i,i+1}=\pi_{i+1} q_{i+1,i}$, $i\in \Nats$. Solving this recursion, we obtain 
\begin{equation}
    \label{eq:invar1}
    \pi_u = \pi_0 \prod_{v=1}^u \frac{d\lambda}{1+v\mu} = (1-\lambda) \prod_{v=1}^u \frac{d\lambda}{1+v\mu},
\end{equation}
where we follow the convention that an empty product is equal to 1. 

Notice that the invariant distribution $\pi$ is a function of $\mu$, which is defined in terms of $\E Q$. But $\E Q$ is the mean of the invariant distribution, $\pi$. Thus, the mean-field equations yield the following fixed-point equation for the mean queue length, $\E Q$
\begin{equation}
    \label{eq:fixed_point}
    \sum_{u=0}^{\infty} \prod_{v=1}^u \frac{d\lambda}{1+v(d-1)\lambda/\E Q}=\frac{1}{1-\lambda}.
\end{equation}
In the theorem below we show that there exists a unique solution to the above fixed point equation for $\lambda \in (0,1)$ and characterise the mean of the resulting distribution.

\begin{theorem}
\label{thm:mean_field}
Consider the mean-field model defined by the transition rates in \eqref{eq:trans_rates}. The following hold:

\begin{enumerate}
    \item For any $\lambda \in (0,1)$, there is a unique $\E Q>0$ which solves~\eqref{eq:fixed_point}. Moreover, if we set $\mu=(d-1)\lambda/\E Q$, then $\pi$ given in \eqref{eq:invar1} is invariant for the Markov process with transition rates given in \eqref{eq:trans_rates}, and the mean of the distribution $\pi$ is equal to $\E Q$. 
    
    \item For any $d\geq 2$ the solution $\E Q$ to~\eqref{eq:fixed_point} satisfies
        $$
        \E Q \leq \frac{-\log(1-\lambda)}{c_d} \quad \forall \; \lambda \in 
        \Bigl( \frac{4d-2}{4d-1},1 \Bigr),
        $$
    where $c_d= \frac{2d-3}{2(d-1)} \log \frac{4d}{4d-1}$.
\end{enumerate}
\end{theorem}


\begin{proof}
To prove the first part, we define for $x>0$ and $u\in \Nats$ the following function 
\begin{equation*}
    \pi_u(x) = (1-\lambda) \prod_{v=1}^u \frac{d\lambda}{1+v(d-1)\lambda x}, \quad \phi(x)=\sum_{u=0}^{\infty} \pi_u(x).
\end{equation*}
Clearly, $x\mapsto \pi_u(x)$ is a strictly decreasing function for each $u\geq 1$, and hence, so is $x\mapsto \phi(x)$. It is also not hard to see that $\phi$ is continuous on $(0,\infty)$, and that, by monotone convergence, $\phi(0)=\lim_{x\downarrow 0} \phi(x)$. Now observe that 
$$
\pi_u(0) =(1-\lambda)(d\lambda)^u, \quad \phi(0)=
\begin{cases} 
(1-\lambda)/(1-d\lambda), & d\lambda<1, \\
+\infty, & d\lambda \geq 1.
\end{cases}
$$
Moreover, for all $x\geq 1$, $\pi_u(x)$ is dominated by $\frac{1-\lambda}{u!}(\frac{d}{(d-1)})^u$, which is a summable 
sequence for all $d\geq 2$. Since $\pi_u(x)$ tends to zero for each $u\geq 1$ as $x$ tends to infinity, while $\pi_0(x)\equiv 1-\lambda$, it follows by dominated convergence that $\phi(x)$ tends to $1-\lambda$. 

Since $\phi(\cdot)$ is a strictly decreasing and continuous function whose limit as $x$ tends to zero is bigger than 1 and whose limit as $x$ tends to infinity is smaller than 1, the equation $\phi(x)=1$ has a unique solution. We denote the solution by $\xi(\lambda,d)$ to make explicit its dependence on the system parameters, $\lambda$ and $d$. Finally, the mean queue length is given by $\E Q(\lambda,d) = 1/\xi(\lambda,d)$. 

To prove the second part of the theorem, observe from~\eqref{eq:invar1} that 
\begin{equation} \label{eq:invar_lowerbd}
\pi_u \geq (1-\lambda) \Bigl( \frac{2d\lambda}{2d-1} \Bigr)^u \quad
\forall \; u\leq u^*=\left\lfloor \frac{(2d-3)\E Q}{2(d-1)\lambda} \right\rfloor, 
\end{equation}
since, for all $v\leq u^*$, we have $1+v\mu \leq d-(1/2)$. It follows that 
$$
    1 \geq \sum_{u=0}^{u^*} \pi_u \geq \pi_{u^*} 
    \geq (1-\lambda) \Bigl( \frac{2d\lambda}{2d-1} \Bigr)^{u^*}.
$$
Taking logarithms on both sides and re-arranging, we get 
$$
u^* \log \frac{2d\lambda}{2d-1} \leq -\log(1-\lambda).
$$
Substituting for $u^*$ from~\eqref{eq:invar_lowerbd}, we see that 
$$
\frac{2d-3}{2(d-1)\lambda} \log \frac{2d\lambda}{2d-1} \E Q \leq -\log(1-\lambda).
$$
Now, for $\lambda \in (\frac{4d-2}{4d-1},1)$, we have $\frac{2d-3}{2(d-1)\lambda} \geq \frac{2d-3}{2(d-1)}$ and 
$\frac{2d\lambda}{2d-1} \geq \frac{4d}{4d-1}$. Hence, it follows from the above that 
$$
c_d \E Q \leq -\log(1-\lambda), \quad \forall \; \lambda \in \Bigl( \frac{4d-2}{4d-1},1 \Bigr).
$$
This establishes the second part of the theorem.
\end{proof}


The above theorem establishes the uniqueness of the distribution
$\pi$ which satisfies the fixed point equation~\eqref{eq:fixed_point} and
also characterises the mean of this distribution.
We conjecture that $\pi$ is indeed the weak limit
of $\pi^n$ as $n \to \infty$. In the next section,
we show numerical evidence to support this conjecture and
throughout the rest of the paper we also make significant progress towards rigorously proving the conjecture.

An important observation to make from Theorem~\ref{thm:mean_field}
is that the mean queue-length under $\pi$ 
(which by Little's law is proportional to
the mean response time of jobs under $\pi$) scales as $O(-\log(1-\lambda))$ as $\lambda \to 1$ for any $d \geq 2$.
This result is to be contrasted with systems with no parallelism, i.e., 
where $d=1$. For $d=1$, the servers behave as independent M/M/1 PS servers. Therefore, the mean queue-length is $\lambda/(1-\lambda)$
which scales as $O(1/(1-\lambda))$ when $\lambda$ approaches $1$.
Thus, we conclude that in heavy traffic, even a small amount of parallelism (e.g., increasing $d$ from $1$ to $2$) can result in a significant reduction
in the mean-response time of jobs.

{\em Transient analysis of $Q()$}: So far we have analysed the stationary distribution $\pi$ of the cavity process $Q()$ which captures the evolution of a single tagged server in the infinite system. Building upon the same ideas, we now obtain the evolution equations for the transient distribution of $Q()$. We conjecture the transient distribution of $Q()$ to be the 
limit of the empirical queue-length process $y^n()$ defined in Section~\ref{sec:model}.
Let $\bar y_k(t)=\prob(Q(t)\geq k)$ for each $k \geq 1$
and $\bar y_0(t)=\prob(Q(t)\geq 0)=1$ for each $t \geq 0$. 
Thus, $\bar y()=(\bar y_k(), k\geq 0)$ denotes the 
complementary cumulative distribution (ccdf) of the process $Q()$
at time $t \geq 0$.
Similarly, for the $n^{\textrm{th}}$ system, we define  $\bar y_k^n(t)=\sum_{l\geq k} y_l^n(t)$ for each $k \geq 1$,
$\bar y_0^n(t)=1$ for each $t \geq 0$, and define the empirical ccdf
of the queue lengths at time $t$ to be $\bar y^n(t)=(\bar y_k^n(t), k\geq 0)$.
To analyse the transient distribution of $Q()$, 
we use the same two assumptions as above, along with the additional assumption
that, at any time $t\geq 0$, the queues at all servers other than tagged server have the same distribution as $Q(t)$.
Then, using the same line of arguments as above, we can write the distribution-dependent transition rates of $Q()$ as
\begin{equation}
    \label{eq:trans_rates_mfe}
    \begin{aligned}
    r_{k,k+1}(\bar y(t))&=d\lambda,\\ 
    r_{k+1,k}(\bar y(t))&=1+(k+1)(d-1)\frac{\bar y_1(t)}{\sum_{l\geq 1} \bar y_l(t)}, 
    \end{aligned}
\end{equation}
for each $k \in \Ints_+$, where for $\bar y_1(t)=0$ the value of $\bar y_1(t)/\sum_{l\geq 1} \bar y_l(t)$ is defined by its limit as $\bar y_1(t) \to 0$. All other transition rates are zero.
Note that these transition rates are consistent with the stationary transition rates defined in~\eqref{eq:trans_rates}.
Thus, using Dynkin's formula, it follows that the ccdf
$\bar y()=(\bar y_k(), k\geq 1)$ of $Q()$ satisfies the following system of ordinary differential equations for each $k \geq 1$:

\begin{equation}
\label{eq:mfe}
    \frac{d\bar{y}_k(t)}{dt}=r_{k,k+1}(\bar y(t))(\bar y_{k-1}(t)-\bar y_k(t))
    -r_{k+1,k}(\bar y(t))
    (\bar y_k(t)-\bar y_{k+1}(t))
\end{equation}
%
We conjecture that the 
stochastic process $\bar y^n()$ converges to $\bar y()$ defined above as $n \to \infty$.

\section{Numerical results} \label{sec:numerical}

In this section, we compare the performance of systems 
with different values of $d$, the number of servers over 
which each job is parallelised, to study the benefits of job 
parallelism. We also compare simulation results with results 
of the mean-field analysis of Section~\ref{sec:mean-field}.
Finally, we provide results for non-exponential service time distributions.

We simulated the system model described in 
Section~\ref{sec:model}, with different numbers of servers, 
$n$, and different values of $\lambda$, the normalised 
arrival rate. We carried out event-driven simulations, 
where each event was taken to be an arrival with probability 
$\lambda/(1+\lambda)$ and a `virtual service' with the 
residual probability $1/(1+\lambda)$. Once the type (arrival 
or virtual service) was allocated to an event, the event was 
then assigned to a server chosen uniformly at random. If the 
event type was an arrival, it then chose another $d-1$ servers 
uniformly at random (without replacement) and joined the queue 
at all $d$ of these servers. If the event type was a virtual 
service, and the queue at the chosen server was empty, nothing 
happened; otherwise, a job was chosen uniformly at random from 
amongst those present at the chosen server, and departed the 
system. A hundred snapshots of the simulated system were taken 
at well-separated epochs to ensure that there was very little 
temporal correlation between snapshots. Queue lengths at all 
$n$ queues were recorded at these times, and used to compute 
empirical cdfs and mean values of queue lengths. 

\begin{figure}[htb]
    \centering
    \includegraphics[width=0.6\textwidth]{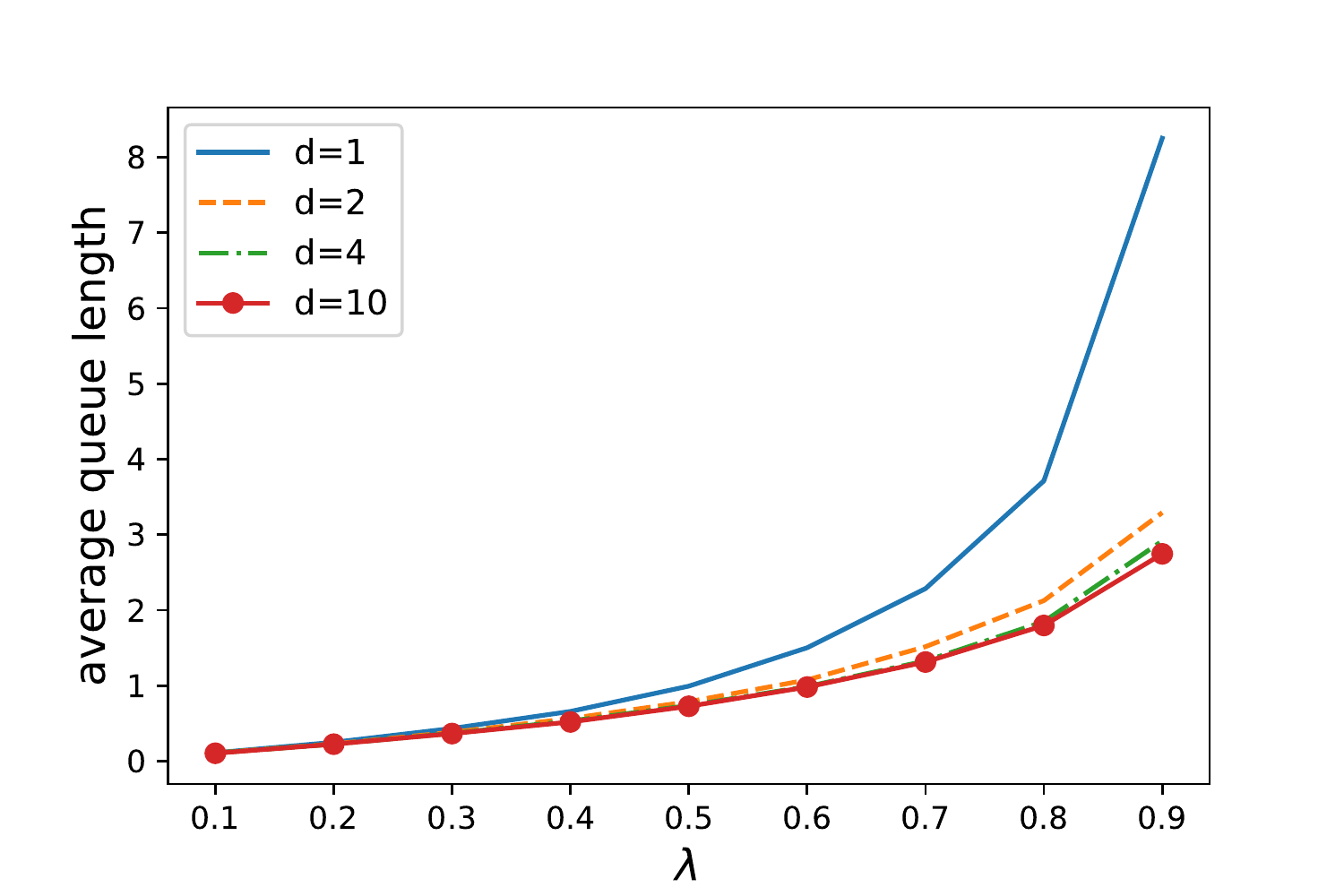}
    \caption{$\E Q$ as a function of $\lambda$}
    \label{fig:mean_queue_length}
\end{figure}

The mean queue length is a key performance measure, 
related to the mean response time via Little's law. 
In Figure~\ref{fig:mean_queue_length}, we plot the mean queue 
length as a function of the normalised arrival rate $\lambda$ 
for a system with $n=100$ servers, for different values of $d$.
Observe that as $d$ increases, the mean queue length 
decreases. The decrease is most significant for higher 
values of $\lambda$ (i.e., as the system approaches 
heavy traffic) and when parallelism increases from 
$d=1$ to $d=2$. This is in accordance with properties 
of the invariant distribution derived in 
Section~\ref{sec:mean-field}. 
Indeed, as shown in Theorem~\ref{thm:mean_field}, for 
$d \geq 2$, the steady-state mean queue length scales 
as $O(\log\brac{\frac{1}{1-\lambda}})$ as $\lambda$ 
increases to $1$. This is a significant decrease 
compared to the $O(1/(1-\lambda))$ scaling for $d=1$.


\begin{table}[htb] 
\renewcommand{\arraystretch}{1.2}
  \centering
      \begin{tabular}{c |c | c | c | c | c |}
    $n$  &  \multicolumn{2}{c|}{$\lambda=0.3$} & \multicolumn{2}{c|}{$\lambda=0.9$}\\   
    \hline
        &  $d=2$ & $d=4$ & $d=2$ & $d=4$\\
    \hline
    $100$ & $0.383$ & $0.368$     & $3.45$ & $2.84$  \\
    $500$ & $0.379$ & $0.366$      & $3.42$ & $2.78$  \\
    Mean-field & $0.379$ &  $0.366$ & $3.34$ & $2.72$\\
    \hline
    \end{tabular}
    \caption{Comparing mean queue lengths from simulations with the mean-field analysis}
    \label{tab:accuracy}
\end{table}

Next, we compare our simulation results to those 
obtained from the mean-field analysis.
In Table~\ref{tab:accuracy} we compare the average queue length obtained from simulations with the average queue length obtained
by solving~\eqref{eq:fixed_point} numerically.
We observe a close match between the results for $n=100$ and $n=500$.

\begin{figure}[htb]
    \centering
    \includegraphics[width=
    0.6\textwidth]{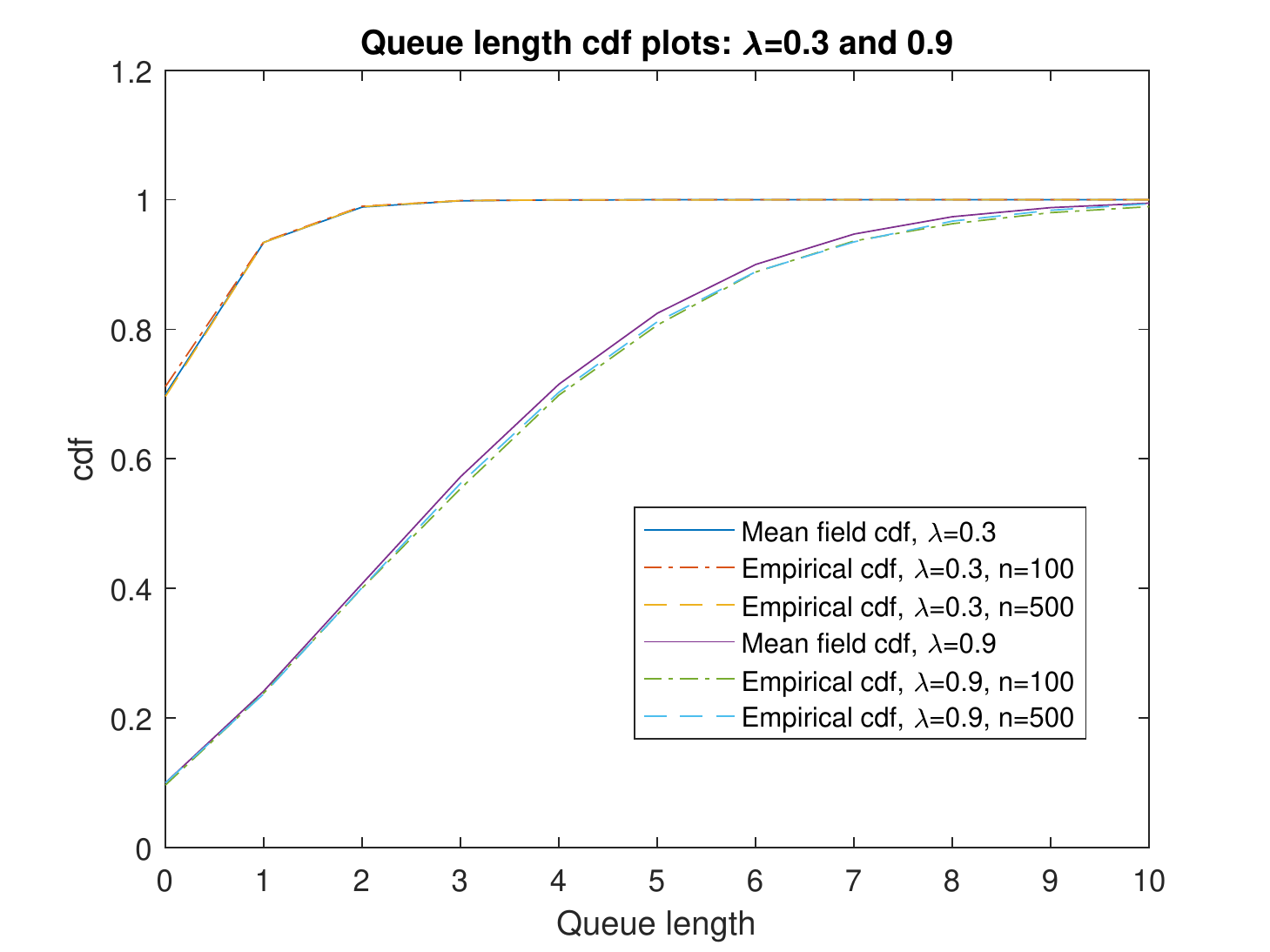}
    \caption{Queue length distribution as a function of $\lambda$ for $d=2$.}
    \label{fig:cdfplots_two_servers}
\end{figure}

\begin{figure}[htb]
    \centering
    \includegraphics[width=0.6\textwidth]{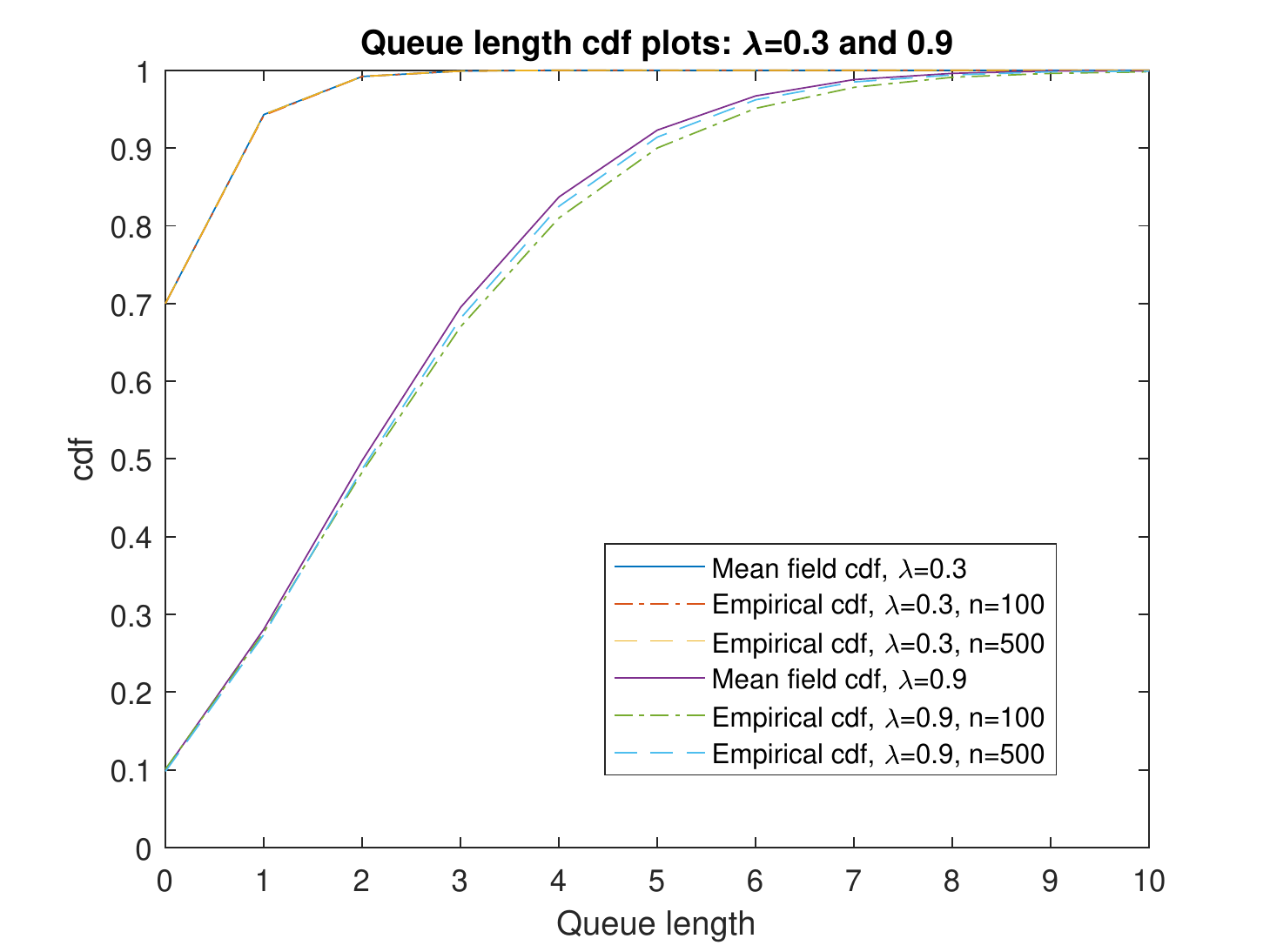}
    \caption{Queue length distribution as a function of $\lambda$ for $d=4$.}
    \label{fig:cdfplots_four_servers}
\end{figure}

Finally, in Figures~\ref{fig:cdfplots_two_servers} and \ref{fig:cdfplots_four_servers}, we plot the empirical 
cdf of the queue length as well as the cdf calculated 
from~\eqref{eq:invar1}, \eqref{eq:fixed_point}. Plots 
are shown for two different arrival rates, $\lambda=0.3$ 
and $\lambda=0.9$, and for the number of servers being 
either $100$ or $500$. The plots in 
Figure~\ref{fig:cdfplots_two_servers} are for $d=2$ (i.e., 
each job is parallelised on 2 servers), while those in 
Figure~\ref{fig:cdfplots_four_servers} are for $d=4$.
For each value of $\lambda$, the figures show that the 
empirical cdfs for $n=100$ and $n=500$ are very close 
to the theoretical cdf from the mean field analysis; 
for $\lambda=0.3$, they are virtually indistinguishable 
in both figures. These results provide additional evidence 
in support of the conjecture that the invariant distribution $\pi^n$
in the system with $n$ servers converges to the solution $\pi$
of the mean-field equations as $n$ tends to infinity.

\begin{figure}[htb]
     \centering
     \begin{subfigure}[b]{0.45\textwidth}
         \centering
         \includegraphics[width=\textwidth]{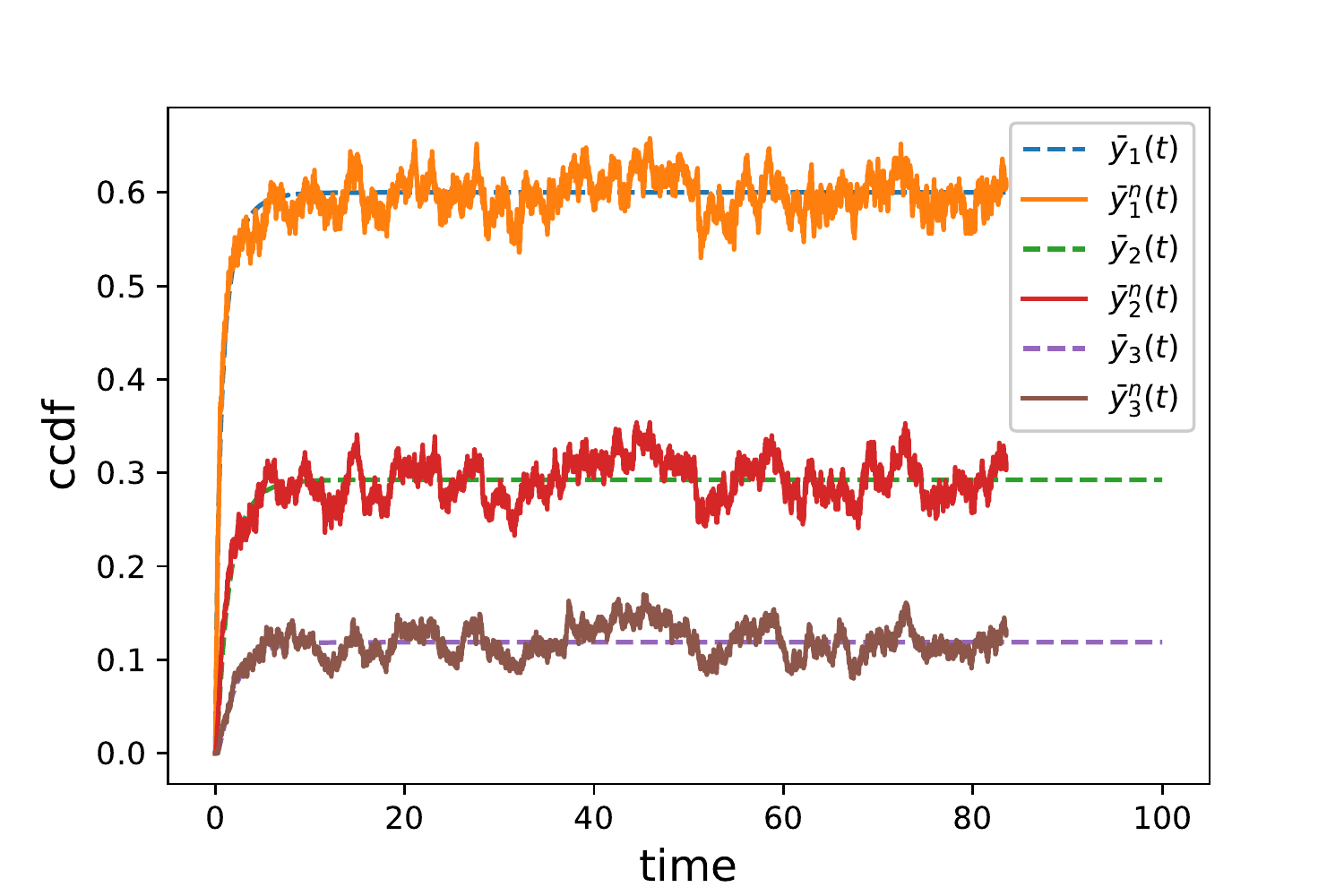}
         \caption{$n=1000$}
         \label{fig:emp_ccdf_n1000}
     \end{subfigure}
     \hfill
     \begin{subfigure}[b]{0.45\textwidth}
         \centering
         \includegraphics[width=\textwidth]{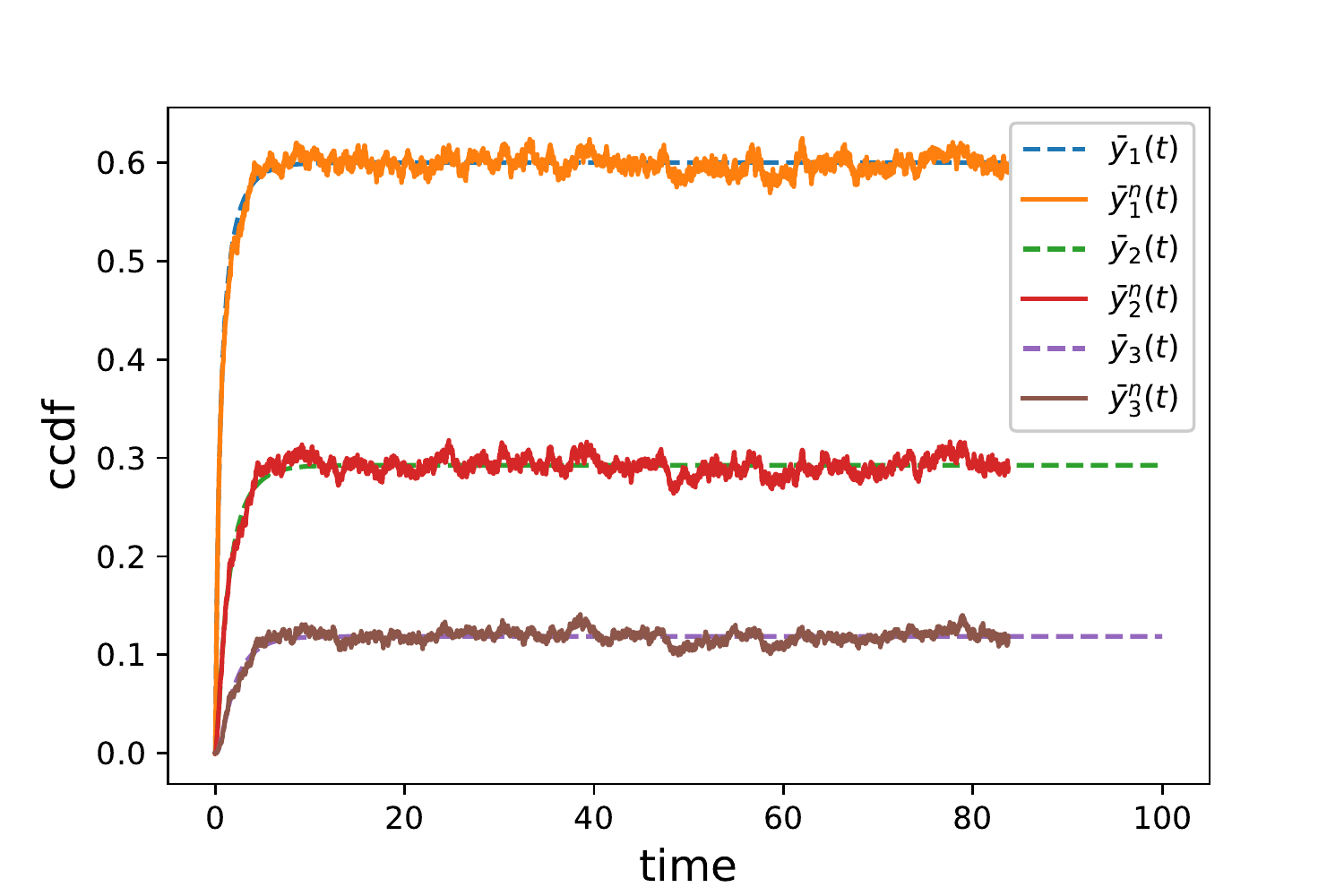}
         \caption{$n=5000$}
         \label{fig:emp_ccdf_n5000}
     \end{subfigure}
     \caption{Comparison of processes $\bar y^n()$ and $\bar y()$ for different values of $n$}
\end{figure}

So far we have presented numerical evidence to support the conjecture
that $\pi^n$ converges to $\pi$ as $n \to \infty$. We now present simulation results to show that similar convergence takes place in the transient regime. More specifically, we show that the process $\bar y^n()=(\bar y^n_k(), k\geq 0)$ converges component-wise to the process $\bar y()=(\bar y_k, k \geq 0)$ defined in Section~\ref{sec:mean-field} as $n \to \infty$. In Figure~\ref{fig:emp_ccdf_n1000} and Figure~\ref{fig:emp_ccdf_n5000} we have compared the first three components of the processes $\bar y^n()$ and $\bar y()$ for $n=1000$ and $n=50000$, respectively. From the figures, it is clear that as $n$ increases the process $\bar y^n()$ concentrates on the mean-field process $\bar y()$. 



\begin{figure}[htb]
    \centering
    \includegraphics[width=0.6\columnwidth]{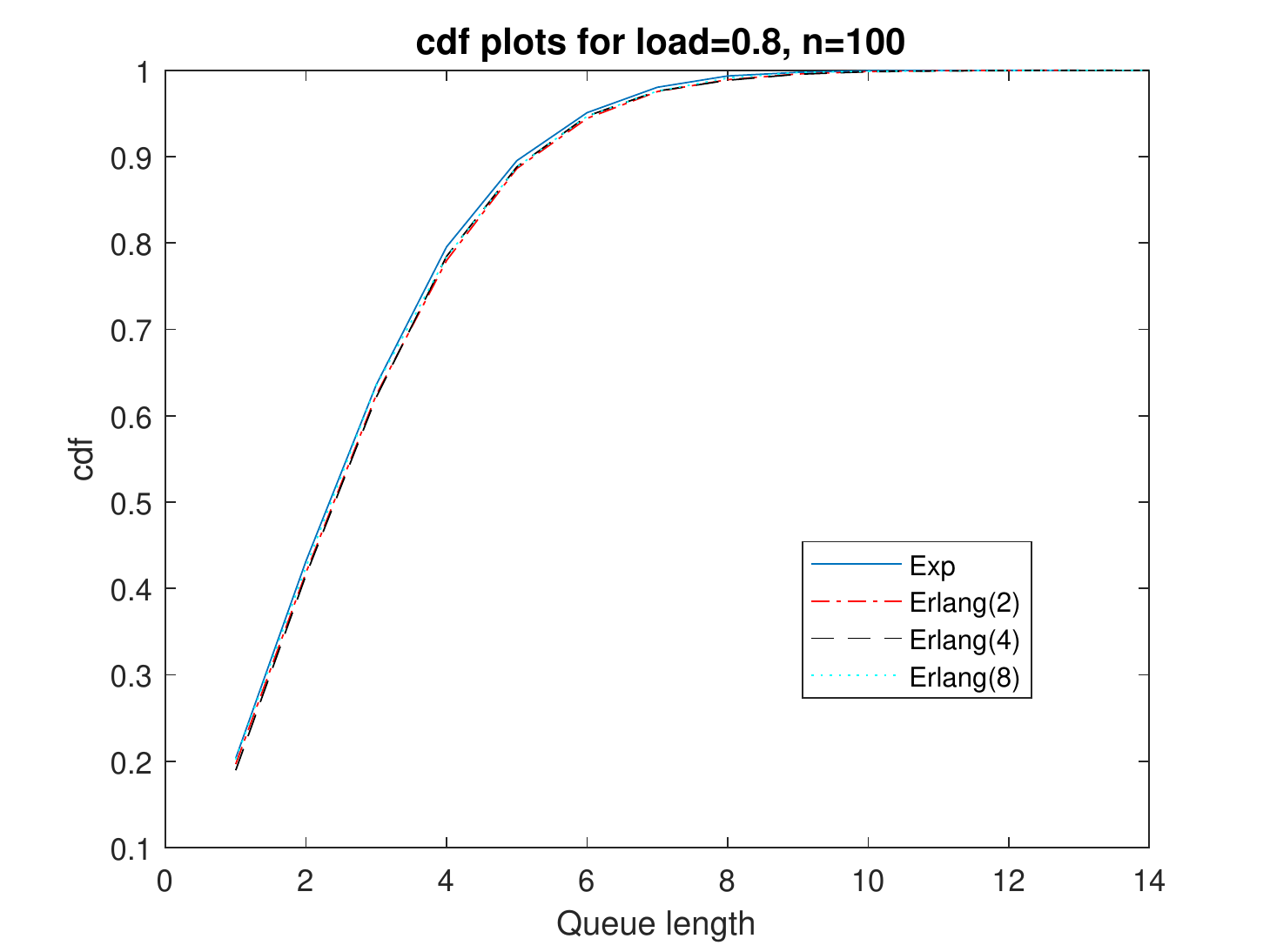}
    \caption{Queue length distribution for Erlang service times, two servers, load=0.8.}
    \label{fig:cdfplot_erlang}
\end{figure}

While the mean-field analysis in the previous section has been 
carried out under the assumption of exponentially distributed 
service time, it is well known that the queue length distribution 
of the $M/G/1-PS$ queue is insensitive to the service time 
distribution but only depends on its mean. This has been proved 
for redundancy models in~\cite{insensitivity_ps_hyperexp}, when 
the service times have hyperexponential distributions with two 
components. It is natural to conjecture that insensitivity holds 
more generally. We now present numerical evidence from simulations 
that this is the case.

\begin{figure}[htb]
    \centering
    \includegraphics[width=0.6\textwidth]{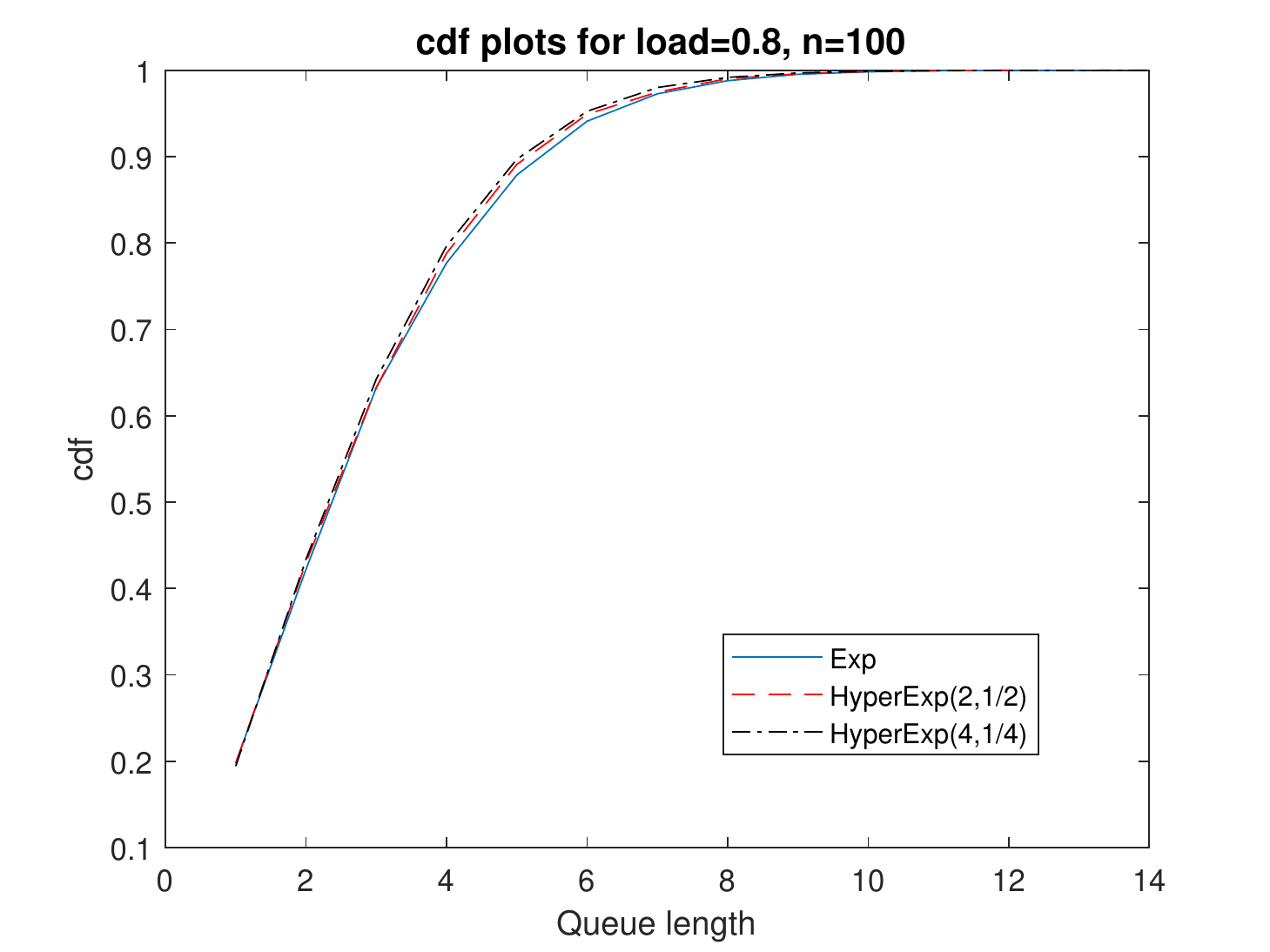}
    \caption{Queue length distribution for hyper-exponential service times, two servers, load=0.8.}
    \label{fig:cdfplot_hyperexp}
\end{figure}

In Figure~\ref{fig:cdfplot_erlang}, we have plotted the empirical cdf 
of the queue length distribution for Erlang service times with 1,2,4 
and 8 phases; if there is 1 phase, service times are exponential. The 
load is kept constant at $0.8$ and each job is assigned to 2 servers. 
The other simulation parameters are similar to the exponential service 
time setting. The empirical cdfs are very close to each other, demonstrating 
that they are insensitive to the service time distribution provided the mean 
is the same. Similar results are seen in Figure~\ref{fig:cdfplot_hyperexp}, 
where we compare empirical queue length cdfs for the Exp(1) service time 
distribution with those for two-component hyperexponentials with parameters 
$(2,1/2)$ and $(4,1/4)$; the weights for the components are chosen to ensure 
the same mean service time of 1. Again, the empirical queue length distribution 
is seen to be insensitive to the service time distribution.

\section{Stability and uniform bounds} \label{sec:stability}

In this section, we show that the system defined in Section~\ref{sec:model} is stable for all $\lambda <1$. 
Moreover, we prove that in the steady state, the expected 
number of jobs at each server is bounded uniformly in $n$.
This implies that the sequence $(\pi^n)_n$ of
steady-state distributions of server occupancy
is tight in $\Ints_+$, which in turn guarantees the
existence of sub-sequential limits. Both these results are essential in showing the . To prove these results,
we analyze the drift of an appropriately defined Lyapunov 
function on the state space of the underlying Markov chain. 

Recall that $X_i^n(t)$ denotes the number of jobs of class $i \in [m]$ in the system at time $t \geq 0$, and that $Q_j^n(t)=\sum_{i \in \partial_j} X_i^n(t)$ denotes the number of ongoing jobs at server $j\in [n]$ at time $t$. 
Note that the dependence of $Q^n()$ on $X^n()$ has not been made explicit in the notation.
The process $X^n()=(X_i^n(),i \in [m])$ is Markov on state space $S^n=\Ints_{+}^{m}$, with transition rates from any state $x \in S^n$
given by

\begin{equation}
\label{eq:transitions}
x \to \begin{cases}
            x+e_i, & \text{at rate  } r^{i,n}_+(x) = \frac{n\lambda}{m},\\
            x-e_i, & \text{at rate  } r^{i,n}_-(x) = x_i \sum_{j \in \partial_i} \frac{1 }{q_j}.
      \end{cases}    
\end{equation}
Here, $q_j=\sum_{i' \in \partial_j}x_{i'}$ is the queue length at server $j$ in state $x$ (whose dependence on $x$ has been suppressed in the notation), while $e_i$ denotes the $m$-dimensional unit vector with the $i^{\textrm{th}}$ component being one. The generator $G_{X^n}$ of the Markov chain $X^n$ is given by

\begin{align}
    \label{eq:gen}
    G_{X^n} f(x) &= \sum_{i \in [m]} \left[r^{i,n}_+(x) \brac{f(x+e_i)-f(x)}  \right.  +\left.  r^{i,n}_-(x) \brac{f(x-e_i)-f(x)}\right],
\end{align}
for $f:S^n \to \Reals$ and $x \in S^n$, where
$r_{\pm}^{i,n}(x)$ are the transition rates from $x$ to $x\pm e_i$
as given by~\eqref{eq:transitions}.
We now establish the following theorem.

\begin{theorem}
\label{thm:stability}
For each $n$, the Markov chain $X^n()$ is positive recurrent
for all $\lambda < 1$. Furthermore, if $Q^n_j(\infty)$, $j\in [n]$, denotes a random variable with the steady-state queue length distribution of server $j$, then we have

\begin{equation*}
    \sup_{n}\E\sbrac{Q^n_j(\infty)} \leq \frac{ \lambda}{1-\lambda}.
\end{equation*}
\end{theorem}

\begin{proof}
We prove the theorem by analyzing the drift of the Lyapunov function $\Psi:S^n \to \Reals_+$ defined as follows
\begin{equation}
\label{eq:psi_def}
    \Psi(x)=\sum_{j \in [n]} (q_j)^2=\sum_{j \in [n]}\Bigl( \sum_{i \in \partial_j} x_i \Bigr)^2.
\end{equation}
%
From~\eqref{eq:gen}, we see that the generator $G_{X^n}$ applied to $\Psi$ and evaluated at $x$ is given by

\begin{align*}
    G_{X^n}\Psi(x)&= \sum_{i \in [m]} \Bigl[ r^{i,n}_+(x)\sum_{j \in \partial_i}(2q_j+1)\Bigr. + \Bigl. r^{i,n}_-(x)\sum_{j \in \partial_i}(-2q_j+1) \Bigr] \\
    &= 2\sum_{i \in [m]} \sum_{j \in \partial_i} [ r^{i,n}_+(x) q_j - r^{i,n}_-(x) q_j]  + d\sum_{i\in [m]} [r^{i,n}_+(x)+ r^{i,n}_-(x)],
\end{align*} 
where we have used the fact that $\sum_{j\in \partial_i} 1=d$ to obtain the second equality.

Note that $\sum_{i\in [m]} r^{i,n}_+(x)$ is just the total arrival rate into the system, which equals $n\lambda$. Similarly, $\sum_{i\in [m]} r^{i,n}_-(x)$ is the total service rate when the system is in state $x$; this is the number of busy servers in state $x$, which we shall denote by $B(x)$. Hence, substituting for the transition rates from~\eqref{eq:transitions}, we obtain from the above that

\begin{equation}
    G_{X^n}\Psi(x)= 2 \frac{n\lambda}{m} \sum_{i \in [m]}\sum_{j \in \partial_i} q_j -2 \sum_{i \in [m]} x_i \sum_{j \in \partial_i} \frac{1}{q_j} \sum_{j \in \partial_i} q_j   + d(n\lambda+B(x)). \label{eq:drift_psi}
\end{equation}   
%
Now, 
\begin{equation}
    \sum_{i \in [m]}\sum_{j \in \partial_i} q_j = \sum_{j\in [n]} \sum_{i\in \partial_j} q_j 
    = \binom{n-1}{d-1} \sum_{j\in [n]} q_j
    = d\binom{n-1}{d-1} \sum_{i\in [m]} x_i, \label{eq:jobcount}
\end{equation} 
where the second equality holds because each server serves $\binom{n-1}{d-1}$ classes; the last equality holds because each job is replicated at $d$ servers, and hence summing the number of jobs at each server gives $d$ times the total number of jobs. 

Next, by the Cauchy-Schwarz inequality,
\begin{equation} \label{eq:bound}
    \sum_{j \in \partial_i} \frac{1}{q_j}\sum_{j \in \partial_i} {q_j} \geq \Bigl( \sum_{j\in \partial_i} 1 \Bigr)^2 = d^2.
\end{equation}
Substituting~\eqref{eq:jobcount} and~\eqref{eq:bound} in~\eqref{eq:drift_psi}, 
we get 

\begin{align}
    G_{X^n}\Psi(x) &\leq 2\lambda d \frac{n}{m} \binom{n-1}{d-1} \sum_{i\in [m]} x_i - 2d^2 \sum_{i\in [m]} x_i +d(n\lambda + B(x)) \nonumber \\
    &= 2d^2 (\lambda-1) \sum_{i\in [m]} x_i +d(n\lambda + B(x)) \nonumber \\
    &= 2d(\lambda-1) \sum_{j\in [n]} q_j + d(n\lambda +B(x)), \label{eq:useful_bound}
\end{align} 
where the first equality follows from the fact that $m=\binom{n}{d}=\frac{n}{d}\binom{n-1}{d-1}$.
%
As the number of busy servers, $B(x)$, is at most $n$, it follows from~\eqref{eq:useful_bound} that, if $\sum_{j \in [n]}q_j > n(1+\lambda)/2(1-\lambda)$, then $G_{X^n} \Psi(X) < 0$. Also note that
\begin{equation}
\label{eq:drift_bound1}
    \sup_{x \in S^n} G_{X^n}\Psi(x) < nd(\lambda+1) < \infty.
\end{equation}
Hence, the drift of $\Psi$ is bounded inside the compact set
$\{x \in S^n: \sum_{j \in [n]} q_j \leq n(1+\lambda)/2(1-\lambda)\}$ and is negative outside this compact set.
Thus, applying the Foster-Lyapunov criterion for stability (see, e.g.,  Proposition D.3 of \cite{Kelly_book})
the process $X^n()$ is positive recurrent.

We now turn to the second statement of the theorem. By Proposition~1 of~\cite{Glynn_bounds}, it follows from \eqref{eq:drift_bound1} that 
\begin{equation}
\label{eq:drift_bound}
    \E[G_{X^n} \Psi(X^n(\infty))] \geq 0.
\end{equation}
Note that we cannot claim that 
$\E[G_{X^n} \Psi(X^n(\infty))] = 0$
as this requires the stronger condition $\E[\Psi(X^n(\infty))] < \infty$ which is not guaranteed by the mere existence of stationary distribution of $X^n()$.
Instead, taking expectation
with respect to the stationary distribution
in~\eqref{eq:useful_bound}, we obtain
\begin{align}
    \E[G_{X^n}\Psi(X^n(\infty))] &\leq -2d(1-\lambda)\E[\sum_{j \in [n]}Q_j^n(\infty)]+nd\lambda+d\E[B(X^n(\infty))] \nonumber \\
    &=-2d(1-\lambda)\E[\sum_{j \in [n]}Q_j^n(\infty)]+2nd\lambda,\label{eq:drift_bound_temp}
\end{align}
where the last equality follows from the fact that, by positive recurrence, the steady-state rate of departures from the system,
$\E[B(X^n(\infty))]$, must be equal to the 
rate of arrivals, $n\lambda$.
Now, substituting the inequality $\E[G_{X^n}\Psi(X^n(\infty))] \geq 0$ in~\eqref{eq:drift_bound_temp}, we get

\begin{align}
    \E[\sum_{j\in [n]}Q_j^n(\infty)] \leq \frac{n\lambda}{1-\lambda},
    \label{eq:q_bound}
\end{align}
from which the second statement of the theorem follows
by noting that for each $j \in [n]$ we have $\E[\sum_{j\in [n]}Q^n_j(\infty)]=n \E[Q^n_j(\infty)]$ due to exchangeablity of the stationary measure.
\end{proof}

\begin{remark}
Theorem~\ref{thm:stability} shows that 
if $\lambda <1$, then stationary queue-length distribution exists and is unique.
Let $\pi^n$ denote the stationary queue-length distribution of an individual server.
From Theorem~\ref{thm:stability} it follows
by Markov inequality that $\sup_n \mathbb{P}(Q^n(\infty) > C) \leq \lambda/(C(1-\lambda))$. Hence, by choosing
$C$ sufficiently large we can guarantee that
$\sup_n \mathbb{P}(Q^n(\infty) > C) \leq \epsilon$
for any $\epsilon > 0$. This shows that the sequence
$(\pi^n)_n$ of individual stationary queue-length
distributions indexed by the system size is tight
in $\mathbb{Z}_+$. Hence, by Prohorov's theorem 
the sequence $(\pi^n)_n$ has subsequential weak limits.
However, to show the convergence of the whole sequence $(\pi^n)_n$
to a unique limit $\pi$, we need to 
further establish that all subsequential limits coincide. 
\end{remark}

\begin{remark}
It follows from~\eqref{eq:q_bound} that 
$$\E[\sum_{i\in [m]}X_i^n(\infty)]\leq  \frac{1}{d}\frac{n\lambda}{1-\lambda}.$$
Combined with Little's law, this implies that
the mean response time of jobs $\bar T_d^n$
in the $n$th system is upper
bounded by
$$\bar T_d^n \leq \frac{1}{d} \frac{1}{1-\lambda}.$$
This bound is achieved with equality for $d=1$. However, for $d\geq 2$ and large $n$, 
the stationary distribution computed using
mean-field heuristics of Section~\ref{sec:mean-field}
suggests that the bound could be lowered
to $O(\frac{1}{\lambda}\log \frac{1}{1-\lambda})$
for $\lambda$ close to $1$.
\end{remark}

\begin{remark}
The upper bound in~\eqref{eq:q_bound} also implies that $$\E[X_i^n(\infty)] \leq \frac{n\lambda}{md(1-\lambda)}.$$
Hence, for $2 \leq d \leq n-1$ we have $\E[X_i^n(\infty)] \to 0$ as $n \to \infty$.
This implies that $X_i^n \to 0$ in probability for each $i \in [m]$.
\end{remark}

\section{Monotonicity} \label{sec:monotone}

In this section, we show that the chain $X^n()$
is monotonic with respect to its starting state.
Specifically, we show that 
if  $X^n()$ and $\bar X^n()$ are two copies
of the same chain with initial states satisfying
$X^n(0) \leq \bar X^n(0)$, then
the sample paths of the chains can be constructed on
the same probability space such that $X^n(t) \leq \bar X^n(t)$
for all $t \geq 0$. Here, $X \leq Y$ for $X,Y \in S^n$ means $X_i \leq Y_i$
for all $i \in [m]$.
The main result of this section is as follows.

\begin{theorem}
\label{prop:coupling}
Consider two Markov chains $X^n()$ and $\bar X^n()$, both
evolving according to the transition rates given by~\eqref{eq:transitions}.
Let $X^n(0) \leq \bar X^n(0)$. Then, there exists a coupling 
between the chains such that
$X^n(t) \leq \bar X^n(t)$ for all $t \geq 0$.
In other words, $X^n(t) \underset{st}{\leq} \bar X^n(t)$ for all $t \geq 0$.
\end{theorem}

\begin{proof}
We call the systems corresponding to
processes $X^n$ and $\bar X^n$ as
smaller and larger systems, respectively. We denote 
the queue lengths in these systems by $Q^n$ and 
$\bar Q^n$ respectively, and observe that if 
$X^n_i(t) \leq \bar X^n_i(t)$ for all $i\in [m]$, then 
$Q^n_j(t) \leq \bar Q^n_j(t)$ for all $j\in [n]$.

The coupling is described as follows:
Let the current instant be denoted by
$t$ and assume that $X^n(t) \leq \bar X^n(t)$.
We shall specify a way of
generating the next event and the time $s$ for the next event
such that the inequality $X^n(s) \leq \bar X^n(s)$
is maintained just after the next 
event has taken place.
An event can be either an arrival or a departure
of a class $i$ job for some $i \in [m]$.

For each class $i \in [m]$, we generate
the time until the next arrival as an exponential
random variable with rate $n\lambda/m$ for both systems.
Hence, for both systems, arrivals of jobs
of each class occur at the same instants.

We now turn to departures. If $X_i(t) < \bar X_i(t)$ for some $i \in [m]$, then we generate the time until the next class $i$ departure as independent exponential random variables with rates $r^{i,n}_-(X^n(t))$ and $r^{i,n}_-(\bar X^n(t))$ in the $X^n$ and $\bar X^n$ systems, respectively. 

Otherwise, if $X_i(t) = \bar X_i(t)$ for some $i \in [m]$, then observe that 
\begin{equation*}
r^{i,n}_-(X^n(t))=X^n_i(t) \sum_{j \in \partial_i} \frac{1}{Q_j^n(t)} \geq \bar X^n_i(t) \sum_{j \in \partial_i} \frac{1}{\bar Q^n_{j}(t)} = r^{i,n}_-(\bar X^n(t)).
\end{equation*}
We first generate
the time till the next class $i$ departure
in the larger system as an exponential
random variable $\bar Z$ with rate $r^{i,n}_-(\bar X^n(t))$.
Then, we generate the time till the next class
$i$ departure for the smaller system as
$Z=\min(\bar Z,Y)$ where $Y$ is another independent
exponential random variable with rate $r^{i,n}_-(X^n(t))-r^{i,n}_-(\bar X^n(t))$, which we 
showed was non-negative. 
(If an exponential random variable has rate 0, we 
formally set it to be $+\infty$.)
%
%
Since $Z \leq \bar Z$, the 
next class $i$ departure
event occurs earlier in the smaller system
than in the larger system.

Once all the event times have been generated
in the way described above, the next event is
taken to be the one whose next event
time is the earliest. Due to the construction
above, it is clear that $X^n(s) \leq \bar X^n(s)$
holds just after the next event time $s$.
This completes the proof of the theorem.
%
\end{proof}

\begin{corollary}
\label{cor:monotone}
If $X^n(0)=0$, then $X^n(s) \underset{st}{\leq} X^n(t)$ for all $0 \leq s \leq t$.
\end{corollary}

\begin{proof}
We construct two processes $Y^n()$ and $\bar Y^n()$ such that
$Y^n(0)=0=X^n(0)$ and $\bar Y^n(0)=X^n(t-s) \geq Y^n(0)$ and both evolve according
to the rates given by \eqref{eq:transitions}.
According to the coupling constructed in Theorem~\ref{prop:coupling} we have 
$X^n(s)=Y^n(s) \underset{st}{\leq} \bar Y^n(s) = X^n(t)$.
\end{proof}

It is important to note that $X^n(t) \leq \bar X^n(t)$
implies that $Q^n(t) \leq \bar Q^n(t)$, where $Q^n$
and $\bar Q^n$ denote the corresponding queue-length processes.
Hence, from Theorem~\ref{prop:coupling} it follows that
$Q^n(t) \leq \bar Q^n(t)$
for all $t \geq 0$ if $X^n(0) \leq \bar X^n(0)$.
Also, from Corollary~\ref{cor:monotone} it follows
that $Q^n(s) \underset{st}{\leq} Q^n(t)$ for all
$0 \leq s \leq t$ if $Q^n(0)=0$.

\section{Uniform convergence to stationary distribution}
\label{sec:unif_conv}

In Section~\ref{sec:stability}, we showed for any fixed 
$n\in \Nats$ that the Markov process $X^n()$ is positive 
recurrent under the condition $\lambda<1$, which we assume 
throughout. Hence, by irreducibility, $X^n()$ has a unique
invariant distribution. Consequently, so do the joint 
and marginal queue length processes $Q^n()$ and $Q^n_j()$. 
Denote by $\pi^n$ the invariant queue length distribution 
at the first server, which is the same as at any other 
server by exchangeability. In the previous section, we 
showed that, if the Markov process $X^n()$ is 
started in the empty state $X^n(0)\equiv 0$, then 
$Q^n_1()$  converges in distribution monotonically 
to $\pi^n$. In this section, we strengthen this result 
by showing that the speed of convergence, made precise 
below, does not depend on $n$. 
Our main result in this section is stated below.

\begin{theorem}
\label{prop:unif_conv}
Let $\pi^{n,0}(t)$ denote the queue-length distribution 
of the first server at time $t\geq 0$, starting from 
the empty system at $t=0$ (i.e., $X^n_i(0)=0$ for all 
$i \in [m]$). Let $\pi^n$ denote the stationary 
queue-length distribution of the first server.
Then, for each $\epsilon >0$, there exists 
$\tau=\tau(\epsilon)$ (not depending on $n$) such that 

\begin{equation*}
    d_{TV}(\pi^{n,0}(t),\pi^n) \leq \epsilon, \text{ for all } t \geq \tau,
\end{equation*}
where $d_{TV}(\cdot,\cdot)$ denotes the total variation
distance.
\end{theorem}
Before providing a detailed proof, we outline the approach. 
We consider two chains $X^n()$ and $\bar X^n()$, both evolving 
according to the transition rates given by~\eqref{eq:transitions}, 
with $X^n()$ being started empty and $\bar X^n()$ being started 
in the invariant distribution. Clearly $\bar X^n(t)$ is in the 
invariant distribution at all times $t$. Denote by $Q^n()$ and 
$\bar Q^n()$ the queue lengths induced by $X^n()$ and $\bar X^n()$, respectively; 
then, for any $t>0$ and $j\in [n]$, $Q^n_j(t)$ has distribution 
$\pi^{n,0}(t)$ while $\bar Q^n_j(t)$ has distribution $\pi^n$.
We couple the chains as described in the last section, so 
that the system started empty is always dominated by the one 
started in equilibrium. Now, 
$$
d_{TV}(\pi^{n,0}(t),\pi_n) \leq \prob(Q^n_1(t)\neq \bar Q^n_1(t)) \leq \E \abs{Q^n_1(t)-\bar Q^n_1(t)},
$$ 
where the first inequality holds for any coupling, and the 
second inequality holds because queue lengths are whole numbers 
and differ by at least one when they differ.
We shall bound the last quantity as a function of $t$. In order 
to do so, we define a distance $W: S^n \times S^n \to \mathbb{R}_+$, 
as follows:

\begin{equation*}
W(x^n,\bar x^n) = \frac{1}{n} \sum_{j=1}^n \abs{\bar q^n_j-q^n_j},
\end{equation*}

where $\bar q^n$ and $q^n$ are the queue lengths corresponding 
to the configurations $\bar x^n$ and $\bar x$ respectively. Note 
that $W$ is not a metric on $S^n \times S^n$ because it may be 
zero for two different configurations which induce the same 
queue lengths. However, it is a metric when restricted to the 
subset $\{ (x^n,\bar x^n)\in S^n\times S^n: x^n \leq \bar x^n \}$, 
which is the subset on which we work.
We have the following lemmas.

\begin{lemma}
\label{lem:coupled_gen}
Consider two Markov chains $X^n()$ and $\bar X^n()$,
both evolving according to the transition rates given by~\eqref{eq:transitions}, with $X^n(0) \leq \bar X^n(0)$. 
Then for any $0 \leq u \leq v$ we have

\begin{equation*}
    \E[W(X^n(v),\bar X^n(v))-W(X^n(u),\bar X^n(u))] 
    = -\frac{d}{n} \sum_{j \in [n]}\int_{u}^v  \prob(\bar Q^n_j(s) >  Q^n_j(s) = 0) ds
\end{equation*}
\end{lemma}

\begin{proof}
We couple the Markov chains $X^n()$ and $\bar X^n()$
as described in Theorem~\ref{prop:coupling}.
Under this coupling we have $\bar Q^n_j(t) \geq Q^n_j(t)$ for all $t \geq 0$. 
Hence, for all $t \geq 0$, we can drop absolute values and write 

\begin{equation*}
    W(X^n(t),\bar X^n(t))=\frac{1}{n}\sum_{j=1}^n (\bar Q^n_j(t) - Q^n_j(t)). 
\end{equation*}
Since $\bar X^n()$ is positive recurrent, we further note that 
$\E[W(X^n(t),\bar X^n(t))] \leq \frac{1}{n} \sum_{j \in [n]} \E[\bar Q^n_j(t)] < \infty$.

Let $G_{X^n,\bar X^n}$
denote the generator of the coupled process $(X^n(),\bar X^n())$. 
This generator can be written in terms
of the generators of the individual processes
as
\begin{align*}
    & G_{X^n,\bar X^n} W(X^n(t),\bar X^n(t))\\
    & =G_{\bar X^n} \frac{1}{n} \sum_{j \in [n]} \bar Q^n_j(t)
    -G_{X^n} \frac{1}{n} \sum_{j \in [n]} Q^n_j(t)\\
    &=-\frac{d}{n}(B(\bar X^n(t))-B(X^n(t))),
\end{align*}
where the last equality follows easily from~\eqref{eq:gen}
by noting
that $G_{X^n} \frac{1}{n} \sum_{j \in [n]} Q^n_j(t)=d\lambda - \frac{d}{n} B(X^n(t))$.
Under the coupling we have $B( X^n) \leq B(\bar X^n) \leq n$.
Hence, $\E [\abs{G_{X^n,\bar X^n} W(X^n(t),\bar X^n(t))}] \leq d$ for all $t \geq 0$. Hence, it follows
that $W(X^n(t),\bar X^n(t))-\int_{0}^t G_{X^n,\bar X^n} W(X^n(s),\bar X^n(s)) ds$ is a martingale
with respect to the natural filtration of 
the coupled Markov chain. $(\bar X^n(), X^n())$.

We can therefore apply Dynkin's formula to conclude the following, for any $0 \leq u \leq v$:
\begin{align*}
    &\E[W(X^n(v),\bar X^n(v))-W(X^n(u),\bar X^n(u))] \\
    &=\int_u^v \E[G_{X^n,\bar X^n} W(X^n(s),\bar X^n(s))] ds \\
    &=-\frac{d}{n}\int_u^v (\E[B(\bar X^n(s))-B(X^n(s))]) ds \\
    &=-\frac{d}{n}\int_u^v \sum_{j \in [n]} (\prob(\bar Q^n_j(s) > 0) - \prob(Q^n_j(s) > 0)) ds \\
    &= -\frac{d}{n} \int_{u}^v \sum_{j \in [n]} \prob(\bar Q^n_j(s) >  Q^n_j(s) = 0) ds,
\end{align*}
where the last equality follows from the fact
that under the coupling described in
Theorem~\ref{prop:coupling} we have
$\prob(Q^n_j(s) >  \bar Q^n_j(s) = 0)=0$.
\end{proof}

\begin{lemma}
\label{lem:bound1}
Consider two Markov chains $X^n()$ and $\bar X^n()$,
both evolving according to the transition rates given by~\eqref{eq:transitions}, with $X^n(0) \leq \bar X^n(0)$. 
Then for each $M_1 > 0$, there exists $\gamma=\gamma(M_1) >0$ and $M_2=M_2(M_1) >0$ such that

\begin{align}
\label{eq:conv_rate}
    \E[W(X^n(u+M_2),\bar X^n(u+M_2))-W(X^n(u),\bar X^n(u))]\nonumber \\ 
    \leq -\frac{\gamma d}{n}\sum_{j \in [n]} \prob(\bar Q^n_j(u) > Q^n_j(u), Q^n_j(u) \leq M_1), 
\end{align}
for all $u \geq 0$.
\end{lemma}

\begin{proof}
Using Lemma~\ref{lem:coupled_gen} we have
\begin{multline*}
  \E[W(X^n(u+M_2),\bar X^n(u+M_2))-W(X^n(u),\bar X^n(u))]\\
  =-\frac{d}{n} \sum_{j \in [n]} \int_{u}^{u+M_2} \prob(\bar Q^n_j(s) >  Q^n_j(s) = 0) ds  
\end{multline*}
Now, for each $s \geq u \geq 0, j \in [n]$
and $M_1 > 0$ we have
\begin{align}
    &\prob(\bar Q^n_j(s) >  Q^n_j(s) = 0) 
    \geq \prob(\bar Q^n_j(u) >  Q^n_j(u), Q^n_j(u) \leq M_1) \nonumber \\ 
    & \hspace{1em} \times   \prob(\bar Q^n_j(s) >  Q^n_j(s) = 0 \vert \bar Q^n_j(u) >  Q^n_j(u), Q^n_j(u) \leq M_1)  \nonumber
\end{align}
Hence, to prove the lemma it suffices to show
that for each $M_1 >0$ and $j \in [n]$,
there exists $M_2=M_2(M_1)$ and $\gamma = \gamma(M_1) >0$ such that
\begin{align*}
    \int_u^{u+M_2} \prob(\bar Q^n_j(s) >  Q^n_j(s) = 0 \vert & \bar Q^n_j(u) >  Q^n_j(u), 
     Q^n_j(u) \leq M_1) ds \geq \gamma.
\end{align*}

We couple the two systems $\bar X^n()$
and $X^n()$ as described in Theorem~\ref{prop:coupling}.
Let $\bar Q^n_j(u) >  Q^n_j(u)$ and $Q^n_j(u) \leq M_1$.
We choose $M_2=4M_1$ and consider the event $A$
over which (i) no arrivals occur at queue $j$
in either of the two systems $\bar X^n()$ and $X^n()$ 
in the interval $[u,u+M_2]$,
(ii) all original jobs present in queue $j$ of the system $X^n()$
at time $u$ leave the system by time $u+M_2/2$, and
(iii) at least one of the original jobs in queue $j$
of the system $\bar X^n$
remains in the system at time $u+M_2$.
The first two events mentioned above
are independent of each other under the coupling
described in Theorem~\ref{prop:coupling} since. The probability of
the first event is $e^{-d\lambda M_2}=e^{-4d\lambda M_1}$.
The probability of the second event
is at least $\frac{1}{2}$. To see this, let $S_i$
be the service time of the $i^{\textrm{th}}$
original job in queue $j$ of the lower system.
The total amount of work done by the $j^{\textrm{th}}$
server in time $M_2/2$ is $M_2/2$.
Hence, at least one of the original
jobs in queue $j$ of the lower system 
remains in the system
after time $M_2/2$ only if
$\sum_{i \in [M_1]} S_i > \frac{M_2}{2}$.
The probability of this event satisfies
$$\prob\brac{\sum_{i \in [M_1]} S_i \geq \frac{M_2}{2}} \leq \frac{2\E[\sum_{i \in [M_1]}S_i]}{M_2}  =\frac{1}{2}.$$
The third event includes the event that 
one of the original jobs in the larger system
at queue $j$ has a service time of at least $d M_2$
which has a probability of $e^{-d M_2}=e^{-4dM_1}$
and is independent of the first two events.
Thus, we have

\begin{align*}
    \int_u^{u+M_2} & \prob(\bar Q^n_j(s) >  Q^n_j(s) = 0 \vert \bar Q^n_j(u) >  Q^n_j(u), \\ 
    & \qquad \qquad \qquad \qquad \qquad Q^n_j(u) \leq M_1) ds  \\
    &\geq \int_{u+M_2/2}^{M_2}  \frac{1}{2} e^{-4dM_1(\lambda+1)} ds \\
    &= M_1 e^{-4dM_1(1+\lambda)}
\end{align*}
This shows that~\eqref{eq:conv_rate} holds
with $\gamma=M_1 e^{-4d M_1(1+\lambda)}$
and $M_2=4M_1$.
\end{proof}


\noindent{\em Proof of Theorem~\ref{prop:unif_conv}}:
We prove Theorem~\ref{prop:unif_conv} using the lemmas stated above.

Let $\bar X^n()$ and $X^n()$ be two Markov chains,
both evolving according to the transition rates
given by~\eqref{eq:trans_rates}.
Let $X^n(0)=0$ and $\bar X^n(0)$ be distributed
according to the stationary distribution
of $X^n()$. We couple the two chains according
to the coupling described in Theorem~\ref{prop:coupling}. By the coupling lemma
we have that
$d_{TV}(\pi_0^n(t),\pi^n) \leq \prob(\bar Q_1^n(t) \neq Q_1(t))$, where $\bar Q^n()$
and $Q^n()$ are the queue length processes
corresponding to two chains $\bar X^n()$
and $X^n()$, respectively.

To prove Theorem~\ref{prop:unif_conv}, we shall show that for each $\epsilon >0$
there exists $\tau(\epsilon)$ such that for
some $u \in (0, \tau(\epsilon)]$ we have
\begin{equation}
  \prob(Q^n_1(u) \neq \bar Q^n_1(u)) \leq \epsilon.  
  \label{eq:suff_cond}
\end{equation}
This will imply the statement of Theorem~\ref{prop:unif_conv} since for any $t \geq \tau(\epsilon)$ we have
\begin{align*}
    \prob(Q^n_j(t) \neq \bar Q^n_j(t))&=\prob(\bar Q_j^n(t)-Q_j^n(t) > 0) \\
    & \leq \prob(\bar Q_j^n(t)-Q_j^n(u) > 0) \\
    & = \prob(\bar Q_j^n(u)-Q_j^n(u) > 0) \\
    & = \prob(Q^n_j(u) \neq \bar Q^n_j(u)) \\
    & \leq \epsilon, \nonumber
\end{align*}
where the first and fourth
lines follow due to the coupling
described in Theorem~\ref{prop:coupling};
the inequality in the second line follows
from Corollary~\ref{cor:monotone} due to
the monotonicity of $Q^n()$; the equality
in the third line follows from the fact
$\bar X^n()$ is a stationary chain.

We now proceed to establish~\eqref{eq:suff_cond}.
We observe the following
\begin{align}
    & \prob(Q_1^n(u) \neq \bar Q_1^n(u)) \nonumber\\
    & = \frac{1}{n}\sum_{j \in [n]}\prob(Q_j^n(u) \neq \bar Q_j^n(u)) \nonumber \\
    & \leq \frac{1}{n}\sum_{j \in [n]}(\prob(Q_j^n(u) \neq \bar Q_j^n(u), Q_j^n(u) \leq M_1) \nonumber\\
    &\hspace{6em} + \prob( Q_j^n(u) \geq M_1))\nonumber \\
    & \leq \frac{1}{n}\sum_{j \in [n]}(\prob(Q_j^n(u) \neq \bar Q_j^n(u), Q_j^n(u) \leq M_1) \nonumber \\
    & \hspace{6em}+ \prob(\bar Q_j^n(u) \geq M_1)) \nonumber \\
    & \leq \frac{1}{n}\sum_{j \in [n]}\prob(Q_j^n(u) \neq \bar Q_j^n(u), Q_j^n(u) \leq M_1) + \frac{\E[\bar Q_j^n(\infty)]}{M_1} \nonumber\\
    & \leq \frac{1}{n}\sum_{j \in [n]}\prob(Q_j^n(u) \neq \bar Q_j^n(u), Q_j^n(u) \leq M_1) + \frac{\lambda}{M_1(1-\lambda)}, \label{eq:ineq_tmp}
\end{align}
where the equality in the first line
follows from the exchangeability 
of the queues in both systems;
the inequality in the third line
follows from Theorem~\ref{prop:coupling};
the inequality in the last line follows
from Theorem~\ref{thm:stability}.

We now show that each term on the RHS of~\eqref{eq:ineq_tmp} can be bounded above $\epsilon/2$ by appropriately choosing $M_1(\epsilon)$ and $\tau(\epsilon)$.
Clearly, the second term on the RHS of~\eqref{eq:ineq_tmp}
can be bounded by $\epsilon/2$ for an appropriate
choice of $M_1=M_1(\epsilon)$. Thus, it remains to  show that the first term is also bounded above by $\epsilon/2$
for this choice of $M_1$ and 
some $\tau(\epsilon)$ and some $u \leq \tau(\epsilon)$. Suppose this is not
true. In particular, assume that
$$
\frac{1}{n}\sum_{j \in [n]}\prob(Q_j^n(u) \neq \bar Q_j^n(u), Q_j^n(u) \leq M_1(\epsilon)) > \frac{\epsilon}{2}
$$
for all $u \in \{0,M_2,2M_2,\ldots,KM_2\}$, where
$M_2=M_2(M_1)$ is chosen according to Lemma~\ref{lem:bound1}
and $K >0$ is some sufficiently large number to be
chosen later.
Then, applying the above inequality in~\eqref{eq:conv_rate} and summing
over all $u \in \{0,M_2,\ldots,KM_2\}$
we obtain
\begin{align*}
    \E \bigl[ & W \bigl( X^n(M_2(K+1)),\bar X^n(M_2(K+1)) \bigr) \\
    &- W \bigl( X^n(0),\bar X^n(0) \bigr) \bigr] 
    \leq -\frac{\gamma d \epsilon(K+1)}{2},
\end{align*}
where $\gamma=\gamma(M_1)$ is as defined
in Lemma~\ref{lem:bound1}.
Thus, we have
\begin{align*}
    &\E[W(X^n(M_2(K+1)),\bar X^n(M_2(K+1)))] \\
    &\leq \E[W(X^n(0),\bar X^n(0))] -\frac{\gamma d \epsilon(K+1)}{2}, \\
    & = \E[\frac{1}{n}\sum_{j \in [n]}\bar Q_j^n(0)]-\frac{\gamma d \epsilon(K+1)}{2}, \\
    & \leq \frac{\lambda}{1-\lambda} -\frac{\gamma d \epsilon(K+1)}{2},
\end{align*}
where the last inequality follows from
Theorem~\ref{thm:stability}.
The RHS of the above is negative for $K = 2\lambda/(\gamma d \epsilon(1-\lambda))$.
This clearly is a contradiction since
by Theorem~\ref{prop:coupling}
we must have $\E[W(X^n(M_2(K+1)),\bar X^n(M_2(K+1)))] \geq 0$. Hence, for some $u \in \{0,M_2,\ldots,KM_2\}$
we must have
$$\frac{1}{n}\sum_{j \in [n]}\prob(Q_j^n(u) \neq \bar Q_j^n(u), Q_j^n(u) \leq M_1(\epsilon)) \leq \frac{\epsilon}{2}$$
if $K = 2\lambda/(\gamma d \epsilon(1-\lambda))$.
Hence, the statement of the theorem follows
by choosing $\tau(\epsilon)=2KM_2=4M_2\lambda/(\gamma d \epsilon(1-\lambda))$. \qed

\section{Static version}
\label{sec:static}

We have so far considered systems in which jobs are replicated at $d$ servers chosen uniformly at random, wherein servers adopt a processor-sharing policy and allocate their effort equally to all customers in their queue. This leaves open the question of whether weighted processor-sharing schemes, either centralised or distributed, could achieve better performance. Motivated by this question, we now consider a static version of the problem.

Consider a system of $n$ servers and $\lfloor \lambda n \rfloor$ jobs, where $\lambda>0$ is a given constant. Each job is replicated across a set of $d$ servers, chosen independently and uniformly at random from all subsets of size $d$. All jobs are of unit size, and all servers work at unit rate. If we identify servers with vertices and jobs with hyperedges, then the allocation of jobs to servers yields a random $d$-uniform hypergraph (each hyperedge contains exactly $d$ vertices) on $n$ vertices, with $\lfloor \lambda n \rfloor$ hyperedges. We now consider the makespan minimisation problem induced by this hypergraph, as defined below.

Consider the $d$-uniform hypergraph (more precisely, a multigraph) with vertex set $V$ and hyperedge set $E$; each $v\in V$ denotes a server and each $e\in E$ denotes a job. Denote by $x(v,e)$ the fraction of its capacity that server $v$ devotes to job $e$. The makespan minimisation problem is the following linear program:
\begin{equation}
    \label{eqn:lp_makespan}
    \begin{aligned}
    \min_{\boldsymbol{x}} \max_{e\in E} \; & \frac{1}{\sum_{v\in E} x(v,e)} \\
    \mbox{subject to } \; & x(v,e) \geq 0 \, \forall \, v\in V, e\in E, \\
    & x(v,e) = 0, \mbox{ if $v\notin e$}, \\
    & \sum_{e:v\in e} x(v,e) \leq 1.
    \end{aligned}
\end{equation}
In the remainder of this section, we derive bounds on the value of this random linear program by relating it to the $k$-core problem on random hypergraphs.

This same model was proposed in~\cite{hajek90}, where the distribution of the load at a typical server under the optimal assignment was studied. The analysis in \cite{hajek90} was non-rigorous, but was made rigorous for the $d=2$ case in~\cite{anantharam16}; they also gave an exact expression for the maximum load, which is the value of the optimisation problem in~\eqref{eqn:lp_makespan}, in terms of a fixed point problem. Here we consider general $d$, but only obtain bounds on the optimum value.

{\bf Definition.} The $k$-core of a hypergraph is the largest vertex induced subgraph in which the degree of each vertex is at least $k$.

The $k$-core can be obtained by recursively deleting all vertices whose degree is smaller than $k$ and the edges which are incident on them. Note that the $k$-core could be empty. As a simple example, the 1-core of a tree is itself, while its 2-core is empty.

\begin{lemma} \label{lem:makespan}
Let $W_{\max}(H)$ denote the value of the linear program \eqref{eqn:lp_makespan} on a $d$-uniform hypergraph $H=(V,E)$. Let $k^*$ denote the largest value of $k$ for which the $k$-core of $H$ is non-empty. Then,
$$
\frac{k^*}{d} \leq W_{\max} \leq k^*.
$$
\end{lemma}

\begin{proof}
The $(k^*+1)$-core of $H$ is empty, i.e., recursively eliminating vertices of degree $k^*$ eliminates all vertices. For a vertex $v$, let $E_v$ denote the set of hyper-edges which are incident on $v$ in the step at which $v$ is eliminated. Then, $|E_v|\leq k^*$. Set $x(v,e)= 1/k^*$ if $e\in E_v$ and set $x(v,e)=0$ otherwise. In other words, assign a fraction $1/k^*$ of the capacity of a vertex to each hyperedge which is incident upon it in the step at which the vertex is eliminated by the core-finding algorithm. Any residual capacity of that vertex is unassigned (wasted). As there at most $k^*$ incident hyperedges when the vertex is eliminated, the capacity constraint, $\sum_{e:v\in e} x(v,e)\leq 1$, is satisfied. Note that this algorithm assigns capacity at least $1/k^*$ to each job or hyperedge since all hyperedges are eventually eliminated. (If multiple vertices contained in a hyperedge are eliminated in the same time step, it receives capacity a multiple of $1/k^*$). Hence, the makespan is at most $k^*$, which proves the upper bound in the theorem.

Next, let $H'$ denote the $k^*$-core of $H$, which is non-empty. Denote by $|V(H')|$ and $|E(H')|$ the total number of vertices and hyperedges in $H'$. We have 
\begin{equation} \label{eq:core_deg}
d|E(H')| = \sum_{v\in H'} deg_{H'} (v) \geq k^* |V(H')|,
\end{equation} 
where $deg_{H'}$ denotes the degree in the subgraph $H'$. Now, the total work corresponding to jobs or hyperedges in $H'$ is $|E(H')|$. As these hyperedges are only incident to vertices in $H'$ (or else they would have been eliminated), they can only be worked on by a total of $|V(H')|$ vertices. Thus, the makespan is at least $|E(H')|/|V(H')|$. The lower bound in the theorem is now immediate from \eqref{eq:core_deg}.
\end{proof}

We now use the following result about $k$-cores of random hypergraphs to obtain asymptotic bounds on the makespan which hold with high probability (w.h.p.), i.e, with probability tending to 1 as $n$ tends to infinity. For $\mu>0$, let $Z_{\mu}$ denote a Poisson random variable with mean $\mu$. Denote by $H(n,p;d)$ the random $d$-uniform hypergraph on $n$ vertices, where each hyperedge on $d$ vertices is present with probability $p$, independent of all others. We have the following:

\begin{theorem} \label{thm:core_threshold}
Fix natural numbers $d$ and $k$ larger than or equal to $2$. Define 
$$
\gamma_k(d) = \inf_{\mu>0} \frac{\mu(d-1)!}{(\prob(Z_{\mu}\geq k-1)^{d-1}},
$$
and notice that it is finite, strictly positive, and increasing in $k$. Suppose $k$ and $d$ are not both equal to 2. Then, $H(n,\mu/n^{d-1}; d)$ contains a non-empty $k$-core w.h.p. if $\mu>\gamma_k$, whereas its $k$-core is empty w.h.p. if $\mu<\gamma_k$.

In the special case $d=k=2$, we have $\gamma_2(2)=1$. It is well known that the random graph $G(n,\mu/n)$ has a non-empty 2-core (i.e., contains a cycle) w.h.p. if $\mu>1$, whereas, for any $\mu>0$, the probability that $G(n,\mu/n)$ contains a cycle is bounded away from zero, uniformly in $n$.
\end{theorem}

The theorem was proved for random graphs (i.e., when $d=2$) by \cite{pittel96} and extended to random uniform hypergraphs by \cite{molloy05}. It was also noted that, by contiguity, the results extend to the random (hyper)graph model parametrised by the number of (hyper)edges rather than their probability, i.e., the $G(n,m)$ and $H(n,m;d)$ models with $m=\binom{n-1}{d-1}p$.

The following result is thus a corollary of Lemma~\ref{lem:makespan} and Theorem~\ref{thm:core_threshold}.

\begin{corollary} \label{cor:makespan}
Fix $d\geq 2$ and consider a system with $n$ servers and $\lambda n$ jobs, where each job is replicated on a subset of $d$ servers chosen uniformly at random and independent of other jobs. If $d=2$ and $\lambda<1/2$, then the makespan is bounded above by 1. Otherwise, if 
$$
\gamma_k(d) < \lambda \cdot d! < \gamma_{k+1}(d),
$$
then, w.h.p., $k$-core is non-empty and the $k+1$-core is empty, and the makespan is bounded between $k/d$ and $k$.
\end{corollary}

\section{Conclusion and Discussions}
\label{sec:discussion}

We have studied an idealised model of job parallelism where each job receives service from $d \geq 2$ processor sharing servers simultaneously. Using a mean-field model, we have studied the average delay experienced by  jobs, in the limit as the number of servers tends to infinity. Our results
show that a significant reduction in the average delay can be obtained near the heavy-traffic limit.
In particular, the average delay scales as $O(\frac{1}{\lambda} \log \frac{1}{1-\lambda})$ as $\lambda \to 1$ for $d\geq 2$. This is a significant reduction compared to the $d=1$ case where the delay is known to be $1/(1-\lambda)$. Numerical results show that the proposed mean-field approximation is accurate even for moderately large system sizes. 

We make significant progress towards rigorously establishing the results obtained from our mean-field model. In particular, we show that the system is stable when the normalised arrival rate $\lambda$ is below the normalised system capacity. We also show that the individual queue lengths are uniformly bounded for all system sizes. This proves the existence of subsequential limits for the sequence 
of stationary queue-length distributions indexed by the system size. In addition, we establish an important monotonicity property required to study the speed of convergence to the stationary distribution. Using the uniform bounds on the queue lengths and the monotonicity property, we show that the rate at which the individual queue length distribution converges to the corresponding stationary distribution does not depend on the system size $n$. 


We have not been able to show 
that the sequence  of stationary queue-length distributions $(\pi^n)_n$ converges as $n$ tends to infinity
to the conjectured distribution $\pi$ derived in Section~\ref{sec:mean-field}. We now outline some 
ideas as to how this might be proved, building upon the results 
that we have established.


The key idea is to bound the distance between
$\pi^n$ and $\pi$ using triangle inequality as follows:
\begin{equation}
    d_{TV}(\pi^n, \pi) \leq d_{TV}(\pi^n, \pi^{n,0}(t))+d_{TV}(\pi^{n,0}(t),\pi^0(t))
    +d_{TV}(\pi^0(t),\pi), \label{eq:triangle_ineq}
\end{equation}
where $\pi^0(t)=(\pi^0_k(t), k\in \Ints_+)$ denotes the distribution of the queue-length process $Q()$ of the tagged server at time $t$ starting from $Q(0)=0$.
Recall from Section~\ref{sec:mean-field} that the transient ccdf of $Q()$ satisfies the mean-field equations~\eqref{eq:mfe}. Hence, we can obtain
$\pi^0(t)$ from the solution $\bar y^0(t)$ of the mean-field equations starting from the initial ccdf $\bar y(0)=(1,0,0,\ldots)$ which corresponds to the empty queue. If this solution exists and is unique, then we can set
$\pi^0_k(t)=\bar y_k(t)- \bar y_{k+1}(t)$ for each $k\geq 0$. We note that by Theorem~\ref{prop:unif_conv} that the first term in the RHS of \eqref{eq:triangle_ineq} can be bounded by any $\epsilon >0$ for sufficiently large $t$ (independent of $n$). Now, if each of the other two terms in the RHS can be bounded by any $\epsilon >0$ for sufficiently large $n$ and $t$, then we can conclude that $d_{TV}(\pi^n,\pi) \to 0$ as $n \to \infty$. But for this to be true we need to show that (i) for each fixed $t$, the distribution $\pi^{n,0}(t)$ converges to $\pi^0(t)$ as $n \to \infty$, i.e., the system started from the empty state converges to the mean-field and (ii) the transient queue-length distribution $\pi^0(t)$ converges to the invarant distribution of the mean-field model $\pi$ as $t \to \infty$. Establishing these two properties seem difficult since the queue-length process $Q^n()$ is non-Markovian for finite $n$. We leave these as future directions to explore.

Finally, there are other interesting directions to explore as well. For example, the case where $d$ varies with $n$ has not been studied in the paper. We believe that our mean-field model could be extended to the case where $d$ increases sufficiently slowly with $n$. Letting $d \to \infty$ in the mean field model
yields an average queue length of $\log(1/(1-\lambda))$. Our numerical experiments suggest that the average queue length is close to this value when $d$ is scaled as $d=\Theta(\log(n))$. It would be interesting prove this rigorously. 
Another possible extension of our model is the case where the underlying graph is not complete. The key challenge in this context would be to find the conditions under which mean-field results still hold.

\bibliographystyle{plain}
\bibliography{interacting_queues}

\begin{thebibliography}{10}

\bibitem{tensorflow_2016}
Mart\'{\i}n Abadi, Paul Barham, Jianmin Chen, Zhifeng Chen, Andy Davis, Jeffrey
  Dean, Matthieu Devin, Sanjay Ghemawat, Geoffrey Irving, Michael Isard,
  Manjunath Kudlur, Josh Levenberg, Rajat Monga, Sherry Moore, Derek~G. Murray,
  Benoit Steiner, Paul Tucker, Vijay Vasudevan, Pete Warden, Martin Wicke, Yuan
  Yu, and Xiaoqiang Zheng.
\newblock Tensorflow: A system for large-scale machine learning.
\newblock In {\em Proc. USENIX Conf. Oper. Sys. Design and Implem.}, OSDI'16,
  page 265–283, 2016.

\bibitem{afanaseva_2020}
Larisa Afanaseva, Elena Bashtova, and Svetlana Grishunina.
\newblock {Stability Analysis of a Multi-server Model with Simultaneous Service
  and a Regenerative Input Flow}.
\newblock {\em Methodology and Computing in Applied Probability},
  22(4):1439--1455, December 2020.

\bibitem{anantharam16}
Venkat Anantharam and Justin Salez.
\newblock {The densest subgraph problem in sparse random graphs}.
\newblock {\em Ann. Appl. Probab.}, 26(1):305 -- 327, 2016.

\bibitem{anton21}
Elene Anton, Urtzi Ayesta, Matthieu Jonckheere, and Ina~Maria Verloop.
\newblock On the stability of redundancy models.
\newblock {\em Operations Research}, 69(5):1540--1565, 2021.

\bibitem{ayesta18}
Urtzi Ayesta, Tejas Bodas, and Ina~Maria Verloop.
\newblock On a unifying product form framework for redundancy models.
\newblock {\em SIGMETRICS Perform. Eval. Rev.}, 46(3):80–81, jan 2019.

\bibitem{bramson2012}
Maury Bramson, Yi~Lu, and Balaji Prabhakar.
\newblock Asymptotic independence of queues under randomized load balancing.
\newblock {\em Queueing Systems}, 71(3):247--292, 2012.

\bibitem{mapreduce}
Jeffrey Dean and Sanjay Ghemawat.
\newblock Mapreduce: Simplified data processing on large clusters.
\newblock {\em Commun. ACM}, 51(1):107–113, jan 2008.

\bibitem{ganesh2012}
Ayalvadi Ganesh, Sarah Lilienthal, D~Manjunath, Alexandre Proutiere, and
  Florian Simatos.
\newblock Load balancing via random local search in closed and open systems.
\newblock {\em Queueing systems}, 71(3):321--345, 2012.

\bibitem{gardner17}
Kristen Gardner, Mor Harchol-Balter, Alan Scheller-Wolf, Mark Velednitsky, and
  Samuel Zbarsky.
\newblock Redundancy-d: The power of d choices for redundancy.
\newblock {\em Oper. Res.}, 65(4):1078–1094, aug 2017.

\bibitem{gardner16}
Kristen Gardner, Samuel Zbarsky, Sherwin Doroudi, Mor Harchol-Balter, Esa
  Hyyti\"{a}, and Alan Scheller-Wolf.
\newblock Queueing with redundant requests: Exact analysis.
\newblock {\em Queueing Syst. Theory Appl.}, 83(3–4):227–259, aug 2016.

\bibitem{Glynn_bounds}
Peter~W Glynn and Assaf Zeevi.
\newblock Bounding stationary expectations of markov processes.
\newblock In {\em Markov processes and related topics: a Festschrift for Thomas
  G. Kurtz}, pages 195--214. Institute of Mathematical Statistics, 2008.

\bibitem{hajek90}
B.~Hajek.
\newblock Performance of global load balancing by local adjustment.
\newblock {\em IEEE Trans. Info. Theory}, 36(6):1398--1414, 1990.

\bibitem{leastloaded2018}
Tim Hellemans and Benny Van~Houdt.
\newblock On the power-of-d-choices with least loaded server selection.
\newblock {\em SIGMETRICS Perform. Eval. Rev.}, 46(1):114, jun 2018.

\bibitem{least_loaded_heavy_traffic}
Tim Hellemans and Benny Van~Houdt.
\newblock Mean waiting time in large-scale and critically loaded power of d
  load balancing systems.
\newblock {\em Proc. ACM Meas. Anal. Comput. Syst.}, 5(2), jun 2021.

\bibitem{joshi_latency_redundant}
Gauri Joshi, Emina Soljanin, and Gregory Wornell.
\newblock Efficient redundancy techniques for latency reduction in cloud
  systems.
\newblock {\em ACM Trans. Model. Perform. Eval. Comput. Syst.}, 2(2), apr 2017.

\bibitem{Kelly_book}
Frank Kelly and Elena Yudovina.
\newblock {\em Stochastic Networks}.
\newblock Cambridge University Press, 2014.

\bibitem{insensitivity_ps_hyperexp}
Grzegorz Kielanski and Benny Van~Houdt.
\newblock On the asymptotic insensitivity of the supermarket model in processor
  sharing systems.
\newblock {\em Proc. ACM Meas. Anal. Comput. Syst.}, 5(2), jun 2021.

\bibitem{kurtz_1970}
Thomas~G. Kurtz.
\newblock Solutions of ordinary differential equations as limits of pure jump
  markov processes.
\newblock {\em Journal of Applied Probability}, 7(1):49–58, 1970.

\bibitem{lee2017}
Kangwook Lee, Nihar~B Shah, Longbo Huang, and Kannan Ramchandran.
\newblock The mds queue: Analysing the latency performance of erasure codes.
\newblock {\em IEEE Transactions on Information Theory}, 63(5):2822--2842,
  2017.

\bibitem{mallick2020}
Ankur Mallick, Malhar Chaudhari, Utsav Sheth, Ganesh Palanikumar, and Gauri
  Joshi.
\newblock Rateless codes for near-perfect load balancing in distributed
  matrix-vector multiplication.
\newblock {\em SIGMETRICS Perform. Eval. Rev.}, 48(1):95–96, jul 2020.

\bibitem{mezard_cavity}
Marc M{\'e}zard and Giorgio Parisi.
\newblock {The Bethe lattice spin glass revisited}.
\newblock {\em {The European Physical Journal B: Condensed Matter and Complex
  Systems}}, 20:217--240, 2001.
\newblock 23 pages, 6 figures.

\bibitem{mitzenmacher01}
M.~Mitzenmacher.
\newblock The power of two choices in randomized load balancing.
\newblock {\em IEEE Trans. Parallel and Distrib. Sys.}, 12(10):1094--1104,
  2001.

\bibitem{molloy05}
Michael Molloy.
\newblock Cores in random hypergraphs and boolean formulas.
\newblock {\em Random Struct. Algorithms}, 27(1):124–135, aug 2005.

\bibitem{pittel96}
Boris Pittel, Joel Spencer, and Nicholas Wormald.
\newblock Sudden emergence of a giant k-core in a random graph.
\newblock {\em Journal of Combinatorial Theory. Series B}, 67(1):111--151, May
  1996.

\bibitem{stability_redundancy_fcfs}
Youri Raaijmakers and Sem Borst.
\newblock Achievable stability in redundancy systems.
\newblock {\em Proc. ACM Meas. Anal. Comput. Syst.}, 4(3), nov 2020.

\bibitem{tail_redundancy_ps}
Youri Raaijmakers, Sem Borst, and Onno Boxma.
\newblock Stability and tail behavior of redundancy systems with processor
  sharing.
\newblock {\em Performance Evaluation}, 147:102195, 2021.

\bibitem{rizk_2016}
Amr Rizk, Felix Poloczek, and Florin Ciucu.
\newblock Stochastic bounds in fork---join queueing systems under full and
  partial mapping.
\newblock {\em Queueing Syst.}, 83(3–4):261–291, aug 2016.

\bibitem{rumayanstev_2017}
Alexander Rumyantsev and Evsey Morozov.
\newblock Stability criterion of a multiserver model with simultaneous service.
\newblock {\em Ann. Oper. Res.}, 252(1):29--39, 2017.

\bibitem{shneer2021}
S.~Shneer and A.L. Stolyar.
\newblock Large-scale parallel server system with multi-component jobs.
\newblock {\em Queueing Syst.}, 98:21--48, 2021.

\bibitem{GoogleBorg_2015}
Abhishek Verma, Luis Pedrosa, Madhukar~R. Korupolu, David Oppenheimer, Eric
  Tune, and John Wilkes.
\newblock Large-scale cluster management at {Google} with {Borg}.
\newblock In {\em Proceedings of the European Conference on Computer Systems
  (EuroSys)}, Bordeaux, France, 2015.

\bibitem{dobrushin96}
N.~D. Vvedenskaya, R.~L. Dobrushin, and F.~I. Karepelevich.
\newblock Queueing system with selection of the shortest of two queues: an
  asymptotic approach.
\newblock {\em Problems Inform. Transmission}, 32:15--27, 1996.

\bibitem{weina_fj_2018}
Weina Wang, Mor Harchol-Balter, Haotian Jiang, Alan Scheller-Wolf, and
  R.~Srikant.
\newblock Delay asymptotics and bounds for multi-task parallel jobs.
\newblock {\em SIGMETRICS Perform. Eval. Rev.}, 46(3):2–7, jan 2019.

\bibitem{weina_multiserv_2021}
Weina Wang, Qiaomin Xie, and Mor Harchol-Balter.
\newblock Zero queueing for multi-server jobs.
\newblock In {\em Proc. 2021 ACM SIGMETRICS}, SIGMETRICS '21, page 13–14,
  2021.

\end{thebibliography}

\end{document}